\documentclass[11pt]{amsart}
\usepackage{amsfonts}
\usepackage{amsmath}
\usepackage{amssymb}
\usepackage{color}
\usepackage{graphicx}
\usepackage{xspace}
\usepackage{axodraw}
\usepackage[all]{xy}
\usepackage{hyperref}
\usepackage{bm}
\usepackage{wasysym}
 
 \font \sevenrm=cmr7
 \font \fiverm=cmr5

 \newcommand{\nc}{\newcommand}

\newtheorem{thm}{Theorem}

\newtheorem{cor}[thm]{Corollary}

\newtheorem{lem}[thm]{Lemma}
\newtheorem{prop}[thm]{Proposition}
\newtheorem{defn}{Definition}
\newtheorem{rmk}[thm]{Remark}

\nc{\tblue}[1]{\textcolor{blue}{#1}}
\nc{\xing}[1]{\tblue{\underline{Xing:}#1 }}

\nc{\ignore}[1]{{}}
\nc{\mrm}[1]{{\rm #1}}
\nc{\dirlim}{\displaystyle{\lim_{\longrightarrow}}\,}
\nc{\invlim}{\displaystyle{\lim_{\longleftarrow}}\,}
\nc{\vep}{\varepsilon} \nc{\ep}{\epsilon}
\nc{\sigmat}{\widetilde\sigma}
\nc{\ostar}{\overline{*}}

\def\fleche#1{\mathop{\hbox to #1 mm{\rightarrowfill}}\limits}
\def\gfleche#1{\mathop{\hbox to #1 mm{\leftarrowfill}}\limits}
\def\inj#1{\mathop{\hbox to #1 mm{$\lhook\joinrel$\rightarrowfill}}\limits}
\def\ginj#1{\mathop{\hbox to #1 mm{\leftarrowfill$\joinrel\rhook$}}\limits}
\def\surj#1{\mathop{\hbox to #1 mm{\rightarrowfill\hskip 2pt\llap{$\rightarrow$}}}\limits}
\def\gsurj#1{\mathop{\hbox to #1 mm{\rlap{$\leftarrow$}\hskip 2pt \leftarrowfill}}\limits}

\def\sshu{{\joinrel{\ \scriptscriptstyle\sqcup\hskip -1.5pt\sqcup}\,}}

\nc{\mchar}{\mrm{Char}}
\nc{\Hom}{\mrm{Hom}}
\nc{\id}{\mrm{id}}

\nc{\remark}{\noindent{\bf{Remark:}}}
\nc{\remarks}{\noindent{\bf{Remarks:}}}

 \nc{\delete}[1]{}
 \nc{\grad}[1]{^{({#1})}}
 \nc{\fil}[1]{_{#1}}

\nc{\BA}{{\Bbb A}} \nc{\CC}{{\Bbb C}} \nc{\DD}{{\Bbb D}}
\nc{\EE}{{\Bbb E}} \nc{\FF}{{\Bbb F}} \nc{\GG}{{\Bbb G}}
\nc{\HH}{{\Bbb H}} \nc{\LL}{{\Bbb L}} \nc{\NN}{{\Bbb N}}
\nc{\PP}{{\Bbb P}} \nc{\QQ}{{\Bbb Q}} \nc{\RR}{{\Bbb R}}
\nc{\TT}{{\Bbb T}} \nc{\VV}{{\Bbb V}} \nc{\ZZ}{{\Bbb Z}}
\nc{\Cal}[1]{{\mathcal {#1}}}
\nc{\mop}[1]{\mathop{\hbox {\rm #1}}\nolimits}
\nc{\smop}[1]{\mathop{\hbox {\sevenrm #1}}}
\nc{\ssmop}[1]{\mathop{\hbox {\fiverm #1}}}
\nc{\mopl}[1]{\mathop{\hbox {\rm #1}}\limits}
\nc{\frakg}{{\frak g}}
\nc{\g}[1]{{\frak {#1}}}
\def \restr#1{\mathstrut_{\textstyle |}\raise-8pt\hbox{$\scriptstyle #1$}}
\def \srestr#1{\mathstrut_{\scriptstyle |}\hbox to
  -1.5pt{}\raise-4pt\hbox{$\scriptscriptstyle #1$}}
  \nc{\llbracket}{[\![}
  \nc{\rrbracket}{]\!]}
\nc{\wt}{\widetilde}
\nc{\wh}{\widehat}
\nc{\un}{\hbox{\bf 1}}
\nc{\redtext}[1]{\textcolor{red}{\tt #1}}
\nc{\bluetext}[1]{\textcolor{blue}{#1}}
\nc{\comment}[1]{[[{\tt {#1}}]] }
\nc{\R}{{\mathbb R}}

\newcommand{\ind}{1\hskip -4pt 1}
\def\shift#1{\mathop{#1}^{\leftarrow}\limits}

\def\semi{\mathrel{\times}\kern -6.5pt\joinrel\mathrel{\raise 1.4pt\hbox{${\scriptscriptstyle |}$}}\kern 2pt}
\def\racine{{\scalebox{0.3}{ 
\begin{picture}(12,12)(38,-38)
\SetWidth{0.5} \SetColor{Black} \Vertex(45,-33){5.66}
\end{picture}}}}

 \def\arbrea{\,{\scalebox{0.15}{
  \begin{picture}(8,55) (370,-248)
    \SetWidth{2}
    \SetColor{Black}
    \Line(374,-244)(374,-200)
    \Vertex(374,-197){9}
    \Vertex(375,-245){12}
  \end{picture}
}}\,}

 \def\arbreba{\,{\scalebox{0.15}{
\begin{picture}(8,106) (370,-197)
    \SetWidth{2}
    \SetColor{Black}
    \Line(374,-193)(374,-149)
    \Vertex(374,-146){9}
    \Vertex(375,-194){12}
    \Line(374,-142)(374,-98)
    \Vertex(374,-95){9}
  \end{picture}
}}\,}

 \def\arbrebb{\,{\scalebox{0.15}{
  \begin{picture}(48,48) (349,-255)
    \SetWidth{2}
    \SetColor{Black}
    \Vertex(375,-252){12}
    \Line(376,-250)(395,-215)
    \Line(373,-251)(354,-214)
    \Vertex(353,-211){9}
    \Vertex(395,-213){9}
  \end{picture}
}}}

\def\arbreca{\,{\scalebox{0.15}{
\begin{picture}(8,156) (370,-147)
    \SetWidth{2}
    \SetColor{Black}
    \Line(374,-143)(374,-99)
    \Vertex(374,-96){9}
    \Vertex(375,-144){12}
    \Line(374,-92)(374,-48)
    \Vertex(374,-45){9}
    \Line(374,-42)(374,2)
    \Vertex(374,5){9}
  \end{picture}
}}\,}

\def\arbrecb{\,{\scalebox{0.15}{
\begin{picture}(48,94) (349,-255)
\SetWidth{2}
\SetColor{Black}
\Line(376,-204)(395,-169)
\Line(373,-205)(354,-168)
\Vertex(353,-165){9}
\Vertex(395,-167){9}
\Vertex(374,-205){9}
\Line(374,-246)(374,-209)
\Vertex(374,-252){12}
\end{picture}}}\,}

\def\arbrecc{\,{\scalebox{0.15}{
 \begin{picture}(48,98) (349,-205)
    \SetWidth{2}
    \SetColor{Black}
    \Vertex(375,-202){12}
    \Line(376,-200)(395,-165)
    \Line(373,-201)(354,-164)
    \Vertex(353,-161){9}
    \Vertex(395,-163){9}
    \Line(353,-160)(353,-113)
    \Vertex(353,-111){9}
  \end{picture}
}}\,}

\def\arbrecd{\,{\scalebox{0.15}{
\begin{picture}(48,52) (349,-251)
    \SetWidth{2}
    \SetColor{Black}
    \Vertex(375,-248){12}
    \Line(376,-246)(395,-211)
    \Line(373,-247)(354,-210)
    \Vertex(353,-207){9}
    \Vertex(395,-209){9}
    \Line(375,-247)(375,-206)
    \Vertex(376,-203){9}
  \end{picture}
 }}\,}

\def\arbreda{\,{\scalebox{0.15}{
\begin{picture}(8,204) (370,-99)
    \SetWidth{2}
    \SetColor{Black}
    \Line(374,-95)(374,-51)
    \Vertex(374,-48){9}
    \Vertex(375,-96){12}
    \Line(374,-44)(374,0)
    \Vertex(374,3){9}
    \Line(374,6)(374,50)
    \Vertex(374,53){9}
    \Line(374,53)(374,98)
    \Vertex(374,101){9}
  \end{picture}
}}\,}

\def\arbredb{\,{\scalebox{0.15}{
\begin{picture}(48,135) (349,-255)
    \SetWidth{2}
    \SetColor{Black}
    \Line(376,-163)(395,-128)
    \Line(373,-164)(354,-127)
    \Vertex(353,-124){9}
    \Vertex(395,-126){9}
    \Vertex(374,-164){9}
    \Line(374,-205)(374,-168)
    \Vertex(374,-207){9}
    \Line(374,-248)(374,-211)
    \Vertex(374,-252){12}
  \end{picture}
}}\,}

\def\arbredc{\,{\scalebox{0.15}{
 \begin{picture}(48,150) (349,-205)
    \SetWidth{2}
    \SetColor{Black}
    \Line(376,-148)(395,-113)
    \Line(373,-149)(354,-112)
    \Vertex(353,-109){9}
    \Vertex(395,-111){9}
    \Line(353,-108)(353,-61)
    \Vertex(353,-59){9}
    \Line(374,-200)(374,-153)
    \Vertex(374,-149){9}
    \Vertex(374,-202){12}
  \end{picture}
}}\,}

\def\arbredd{\,{\scalebox{0.15}{
 \begin{picture}(48,99) (349,-251)
    \SetWidth{2}
    \SetColor{Black}
    \Line(376,-199)(395,-164)
    \Line(373,-200)(354,-163)
    \Vertex(353,-160){9}
    \Vertex(395,-162){9}
    \Vertex(376,-156){9}
    \Vertex(376,-248){12}
    \Line(375,-245)(375,-204)
    \Line(375,-200)(375,-159)
    \Vertex(375,-201){9}
  \end{picture}
}}\,}

\def\arbrede{\,{\scalebox{0.15}{
 \begin{picture}(48,153) (349,-150)
    \SetWidth{2}
    \SetColor{Black}
    \Vertex(375,-147){12}
    \Line(376,-145)(395,-110)
    \Line(373,-146)(354,-109)
    \Vertex(353,-106){9}
    \Vertex(395,-108){9}
    \Line(353,-105)(353,-58)
    \Vertex(353,-56){9}
    \Line(353,-52)(353,-5)
    \Vertex(353,-1){9}
  \end{picture}
}}\,}

\def\arbredf{\,{\scalebox{0.15}{
\begin{picture}(48,98) (349,-205)
    \SetWidth{2}
    \SetColor{Black}
    \Vertex(375,-202){12}
    \Line(376,-200)(395,-165)
    \Line(373,-201)(354,-164)
    \Vertex(353,-161){9}
    \Vertex(395,-163){9}
    \Line(353,-160)(353,-113)
    \Vertex(353,-111){9}
    \Line(395,-159)(395,-112)
    \Vertex(395,-111){9}
  \end{picture}
}}\,}

\def\arbredz{\,{\scalebox{0.15}{
  \begin{picture}(68,88) (329,-215)
    \SetWidth{2}
    \SetColor{Black}
    \Vertex(375,-212){12}
    \Line(376,-210)(395,-175)
    \Line(373,-211)(354,-174)
    \Vertex(353,-171){9}
    \Vertex(395,-173){9}
    \Line(351,-168)(332,-131)
    \Line(355,-168)(374,-133)
    \Vertex(333,-131){9}
    \Vertex(374,-131){9}
  \end{picture}
}}\,}

\def\arbredg{\,{\scalebox{0.15}{
\begin{picture}(48,98) (349,-205)
    \SetWidth{2}
    \SetColor{Black}
    \Vertex(375,-202){12}
    \Line(376,-200)(395,-165)
    \Line(373,-201)(354,-164)
    \Vertex(353,-161){9}
    \Vertex(395,-163){9}
    \Line(375,-201)(375,-160)
    \Vertex(376,-157){9}
    \Vertex(376,-111){9}
    \Line(375,-155)(375,-114)
  \end{picture}
}}\,}

\def\arbredh{\,{\scalebox{0.15}{
 \begin{picture}(90,46) (330,-257)
    \SetWidth{2}
    \SetColor{Black}
    \Vertex(375,-254){12}
    \Line(376,-252)(395,-217)
    \Vertex(395,-215){9}
    \Line(374,-254)(335,-226)
    \Vertex(334,-224){9}
    \Line(375,-252)(356,-215)
    \Vertex(355,-215){9}
    \Line(374,-255)(417,-227)
    \Vertex(418,-225){9}
  \end{picture}
}}\,}

\def\arbrefr{\,{\scalebox{0.15}{
   \begin{picture}(160,226) (64,-79)
    \SetWidth{2.0}
    \SetColor{Black}
    \Vertex(96,34){16}
    \Vertex(192,34){16}
    \Vertex(144,-62){16}
    \Vertex(192,130){16}
    \Line(96,34)(144,-62)
    \Line(192,34)(144,-62)
    \Line(192,130)(192,50)
     \SetWidth{1.0}
    \Line(144,-62)(144,18)
    \Line(96,34)(64,98)
    \Line(128,98)(96,34)
    \Line(192,34)(224,98)
  \end{picture}}}\,}

\def\arbrefrpetit{\,{\scalebox{0.10}{
   \begin{picture}(160,226) (64,-79)
    \SetWidth{2.0}
    \SetColor{Black}
    \Vertex(96,34){16}
    \Vertex(192,34){16}
    \Vertex(144,-62){16}
    \Vertex(192,130){16}
    \Line(96,34)(144,-62)
    \Line(192,34)(144,-62)
    \Line(192,130)(192,50)
     \SetWidth{1.0}
    \Line(144,-62)(144,18)
    \Line(96,34)(64,98)
    \Line(128,98)(96,34)
    \Line(192,34)(224,98)
  \end{picture}}}\,}
  
  \def\arbrefrpetitb{\,{\scalebox{0.10}{
   \begin{picture}(160,226) (64,-79)
    \SetWidth{2.0}
    \SetColor{Black}
    \Vertex(96,34){16}
    \Vertex(192,34){16}
    \Vertex(144,-62){16}
    \Vertex(192,130){16}
    \Line(96,34)(144,-62)
    \Line(192,34)(144,-62)
    \Line(192,130)(192,50)
     \SetWidth{1.0}
    \Line(144,-62)(144,18)
    \Line(144,-62)(48,18)
    \Line(128,98)(96,34)
    \Line(192,34)(224,98)
  \end{picture}}}\,}
  
  \def\arbrefrpetitc{\,{\scalebox{0.10}{
   \begin{picture}(160,226) (64,-79)
    \SetWidth{2.0}
    \SetColor{Black}
    \Vertex(96,34){16}
    \Vertex(192,34){16}
    \Vertex(144,-62){16}
    \Vertex(192,130){16}
    \Line(96,34)(144,-62)
    \Line(192,34)(144,-62)
    \Line(192,130)(192,50)
     \SetWidth{1.0}
    \Line(144,-62)(144,18)
    \Line(128,98)(96,34)
    \Line(192,34)(224,98)
    \Line(192,34)(160,98)
  \end{picture}}}\,}
  
   \def\arbrefrpetitd{\,{\scalebox{0.10}{
   \begin{picture}(160,226) (64,-79)
    \SetWidth{2.0}
    \SetColor{Black}
    \Vertex(96,34){16}
    \Vertex(192,34){16}
    \Vertex(144,-62){16}
    \Vertex(192,130){16}
    \Line(96,34)(144,-62)
    \Line(192,34)(144,-62)
    \Line(192,130)(192,50)
     \SetWidth{1.0}
    \Line(144,-62)(144,18)
    \Line(128,98)(96,34)
    \Line(192,34)(224,98)
    \Line(192,210)(192,130)
  \end{picture}}}\,}

  \def\arbrefrun{\,{\scalebox{0.10}{
   \begin{picture}(160,226) (64,-79)
    \SetWidth{2.0}
    \SetColor{Black}
    \Vertex(96,34){16}
    \Vertex(192,34){16}
    \Vertex(144,-62){16}
    \Vertex(192,130){16}
    \Line(96,34)(144,-62)
    \Line(192,34)(144,-62)
    \Line(192,130)(192,50)
     \SetWidth{1.0}
    \Line(144,-62)(144,18)
    \Line(128,98)(96,34)
    \Line(192,34)(224,98)
  \end{picture}}}\,}

  \def\arbrefrdeux{\,{\scalebox{0.10}{
   \begin{picture}(160,226) (64,-79)
    \SetWidth{2.0}
    \SetColor{Black}
    \Vertex(96,34){16}
    \Vertex(192,34){16}
    \Vertex(144,-62){16}
    \Vertex(192,130){16}
    \Line(96,34)(144,-62)
    \Line(192,34)(144,-62)
    \Line(192,130)(192,50)
     \SetWidth{1.0}
    \Line(96,34)(64,98)
    \Line(128,98)(96,34)
    \Line(192,34)(224,98)
  \end{picture}}}\,}

\def\arbrefrtrois{\,{\scalebox{0.10}{
   \begin{picture}(160,226) (64,-79)
    \SetWidth{2.0}
    \SetColor{Black}
    \Vertex(96,34){16}
    \Vertex(192,34){16}
    \Vertex(144,-62){16}
    \Vertex(192,130){16}
    \Line(96,34)(144,-62)
    \Line(192,34)(144,-62)
    \Line(192,130)(192,50)
     \SetWidth{1.0}
    \Line(144,-62)(144,18)
    \Line(96,34)(64,98)
    \Line(128,98)(96,34)
  \end{picture}}}\,}

\begin{document}
\title[Hopf algebra of Multi-indices]{Free Novikov algebras and the Hopf algebra of decorated multi-indices}

\author{Zhicheng Zhu}
\address{School of Mathematics and Statistics, Lanzhou University
Lanzhou, 730000, China}
\email{zhuzhch16@lzu.edu.cn}

\author{Xing Gao}
\address{School of Mathematics and Statistics, Lanzhou University
Lanzhou, 730000, China;
School of Mathematics and Statistics
Qinghai Nationalities University, Xining, 810007, China;
Gansu Provincial Research Center for Basic Disciplines of Mathematics
and Statistics, Lanzhou, 730070, China
}
\email{gaoxing@lzu.edu.cn}

\author{Dominique Manchon}
\address{Laboratoire de Math\'{e}matiques Blaise Pascal, CNRS-Universit\'{e} Clermont-Auvergne, 3 place Vasar\'ely, CS 60026, 63178 Aubi\`ere, France}
\email{Dominique.Manchon@uca.fr}
 
\begin{abstract}
We propose a combinatorial formula for the coproduct in a Hopf algebra of decorated multi-indices which recently appeared in the literature, which can be briefly described as the graded dual of the enveloping algebra of the free Novikov algebra generated by the set of decorations. Similarly to what happens for the Hopf algebra of rooted forests, the formula can be written in terms of admissible cuts. We also prove a combinatorial formula for the extraction-contraction coproduct for undecorated multi-indices, in terms of a suitable notion of covering subforest.
\end{abstract}
\date{April 11th 2024}
\maketitle 
{\bf Keywords:} Multi-indices, Novikov algebra, Hopf algebra, symmetry factor, pre-Lie algebra.

{\bf MSC classification:} 05C05, 16T30, 17A30.
\section{Introduction}\label{sect:intro}
A Hopf algebra of multi-indices appeared in the recent work of P. Linares, F. Otto and M. Tempelmayr on regularity structures \cite{LOT23}, as a new combinatorial tool for handling rough partial differential equations. This new approach was continued in \cite{LO22, BL2023, BD2023}, and adapted to the framework of rough paths and rough differential equations in \cite{L23, BEH2024}. We will stick to this somewhat simpler framework, which permits us to look at $A$-decorated multi-indices: contrarily to what happens in general with regularity structures, the set $A$ of decorations does not play an active role in the definition of the coproduct.\\

Recalling that the free pre-Lie algebra generated by $A$ is the linear span of $A$-decorated rooted trees \cite{CL01, DL02}, the canonical surjective map $\Phi$ from the free pre-Lie algebra $\mop{PL}(A)$ generated by $A$ onto the free Novikov algebra $\mop{N}(A)$ yields, by transposition and multiplicative extension, a canonical embedding $\jmath$ of the multi-index Hopf algebra $\mathcal H_{\smop{LOT}}^A$ into the Hopf algebra $\mathcal H_{\smop{BCK}}^A$ of $A$-decorated rooted forests\footnote{In the undecorated case $A=\{*\}$, the Hopf algebra $\mathcal H^A_{\ssmop{LOT}}$  is a quotient of the original Hopf algebra of \cite{LOT23}. A family of generators is discarded due to the fact that we do not deal with regularity structures.}. This conceptually simple algebraic fact allows a grouping of the terms in the tree expansion \cite{Gubi2010} of the solution of a rough differential equation driven by a H\"older continuous path $X:[0;T]\to \mathbb R^d$ corresponding to any choice of a rough path or branched rough path $\mathbb X$ over $X$ \cite{L98, LCL04, LV07, Gubi2010}. The grouping is performed according to multi-indices, considering for each multi-index $M$ the set of trees $t$ such that $\Phi(t)=M$. Along these lines, a multi-index rough path is a two-parameter family of characters of $\mathcal H_{\smop{LOT}}^A$ subject to Chen's lemma and suitable estimates \cite{LOT23, M24}. \\

The above description as a graded dual, although perfectly rigorous, is not completely explicit, as it depends on a choice of pairing on the symmetric algebra\footnote{The enveloping algebra is identified with the symmetric algebra by the Guin-Oudom construction \cite{GO08} associated to the pre-Lie structure.} of $\mop N(A)$. We propose here a pairing, carefully chosen in order to take symmetry factors into account, which yields an explicit formula for the coproduct. The formula \eqref{coprod-N-bis} thus obtained, reminiscent to the Connes-Kreimer formula for rooted forests, involves a suitable notion of admissible cut for multi-indices.\\

The paper is organised as follows: after a quick reminder on combinatorial Hopf algebras (Paragraph \ref{par:cha}), we recall the construction of the rooted forest Hopf algebra in Paragraph \ref{par:bck}. The multi-index Hopf algebra is introduced in Paragraph \ref{par:nov}. We introduce the pairing in Paragraph \ref{par:pair} and the corresponding Hopf algebra embedding in Section \ref{sect:pair}. Section \ref{sect:bc} introduces graftings and cuts for trees and forests with free edges. The main result of this section is Proposition \ref{sym-cut-graft}, which carefully counts cuts and graftings with the help of symmetry factors. We prove a recursive formula for the coproduct in Paragraph \ref{par:rec}, and finally state and prove our main result (Theorem \ref{explicit} and Equation \eqref{coprod-N-bis}) in the last Paragraph \ref{par:main}. The last section is devoted to a combinatorial formula for the extraction-contraction coproduct (Theorem \ref{ec-main}).\\

Our work is close in spirit to \cite{JZ2023}, where an algebraic formula for the coproduct for the full multi-index Hopf algebra of \cite{LOT23} is given in terms of a different pairing: see (2.16) and (3.12) therein. Our choice of pairing yields a completely combinatorial interpretation of the coproduct of the multi-index Hopf algebra $\mathcal H_{\smop{LOT}}^A$. It would be interesting to see whether such a combinatorial interpretation for the regularity structure Hopf algebra of  \cite{LOT23} and \cite{JZ2023} is available.\\

Let us finally mention the recent work by Y. Bruned and Y. Hou \cite{BH24}, in which the authors give explicit formulae for both coproducts starting from explicit expressions for both corresponding Grossman-Larson products and dualizing them. Their choice of pairing uses a different notion of symmetry factor for a multi-index. Although both approaches give equivalent results, different coefficients therefore appear in the explicit expressions. Our combinatorial interpretation of both coproducts in terms of admissible cuts, and respectively covering subforests, however relies on our convention for symmetry factors in an essential way.\\

\noindent\textbf{Acknowledgements:} We thank Lorenzo Zambotti for his careful reading and his most helpful comments on a previous version of this work. We also thank Yvain Bruned, Pierre Catoire, Lo\"\i c Foissy, Yingtong Hou, Jean-David Jacques and Pablo Linares for illuminating discussions. The first author is supported by the China Scholarship Council (File No. 202306180070). The second author is supported by the National Natural Science Foundation of China (12071191) and Innovative Fundamental Research Group Project of Gansu Province (23JRRA684). Work partially supported by Agence Nationale de la Recherche, ANR-20-CE40-0007 Combinatoire Alg\'ebrique, R\'esurgence, Probabilit\'es Libres et Op\'erades. 
\section{The Hopf algebras of decorated multi-indices}
\subsection{Reminder on combinatorial Hopf algebras}\label{par:cha}
There seems to be no consensus on what a combinatorial Hopf algebra is, despite the vitality and fruitfulness of this research topic. The following \textsl{ad hoc} definition comes from \cite{CEMM20}. A combinatorial Hopf algebra is
\begin{itemize}
\item A graded connected Hopf algebra on a field $\mathbb K$ containing the rationals
\[\mathcal H=\bigoplus_{n\ge 0} \mathcal H_n, \hskip 12mm \mathcal H_0=\mathbb K.\mathbf 1,\]
\item together with a homogeneous linear basis $\mathcal B$ such that the structure constants $\kappa_{ab}^{c}$ and $\kappa_a^{bc}$ of both multiplication and comultiplication are non-negative integers:
\[ab=\sum_{c\in\mathcal B} \kappa_{ab} ^c c,\hskip 12mm \Delta(a)=\sum_{b,c\in \mathcal B}\kappa_a^{bc}b\otimes c.\]
\end{itemize}
Some more requirements can be added, such as
\begin{itemize}
\item Moderate growth condition, e.g.
\[\hbox{dim}\mathcal H_n\le CK^n \hbox{ for some } C,K>0.\]
\item Nondegeneracy condition:
\begin{equation}\label{nondeg}
\mathcal B\cap\hbox{Prim}(\mathcal H)=\mathcal B\cap \mathcal H_1,
\end{equation}
where $\hbox{Prim}(\mathcal H)$ is the set of primitive elements. Recall that an element $x\in \mathcal H$ is primitive if and only if $\Delta(x)=x\otimes \mathbf 1+\mathbf 1\otimes x$.
\end{itemize}
A morphism between two combinatorial Hopf algebras $(\mathcal H,\mathcal B)$ and $(\mathcal H',\mathcal B')$ is a graded Hopf algebra morphism $\varphi:\mathcal H\to\mathcal H'$ such that, for any $b\in \mathcal B$, the image $\varphi(\mathcal B)$ is a finite linear combination of elements of $\mathcal B'$ with nonnegative integer coefficients. This organises those combinatorial Hopf algebras into a category.
\ignore{
\begin{defn}\cite{CEMM20}
Let $(\mathcal H,\mathcal B)$ be a nondegenerate combinatorial Hopf algebra in the sense of the definition above. The {\bf inverse factorial character} $q:\mathcal H\to \mathbf k$ associated to these data is recursively defined by $q(b)=1$ for any $b\in \mathcal B_1=\mathcal B\cap\mathcal H_1$ and, for any homogeneous $x$ of degree $\ge 2$,
\[q(x)=\frac{1}{2^{|x|}-2}\sum_{(x)}q(x')q(x'').\]
\end{defn}
\noindent In particular, we have in the shuffle Hopf algebra of words:
\[q(w)=\frac{1}{|w|!},\]
which justifies the terminology. Notice that $q(b)>0$ holds for any $b\in\mathcal B$ provided the nondegeneracy condition \eqref{nondeg} is fulfilled.
}
\subsection{Pre-Lie algebras and the Butcher--Connes--Kreimer Hopf algebra}\label{par:bck}
A left pre-Lie algebra (here, over a field $\mathbb K$ containing the rationals) is a $\mathbb K$-vector space $P$ together with a bilinear product $\rhd$ such that
\begin{equation}\label{pl}
x\rhd(y\rhd z)-(x\rhd y)\rhd z=y\rhd(x\rhd z)-(y\rhd x)\rhd z
\end{equation}
holds for any $x,y,z\in P$. Pre-Lie algebras, which can be traced back to
the work of A. Cayley \cite{Cay}, are sometimes called {\sl Vinberg algebras\/}, as they appeared explicitly for the first time
in the work of E. B. Vinberg \cite{Vinberg63} under the name \textsl{left-symmetric
algebras} on the classification of
homogeneous cones. They appeared independently at the same time in the work of
M. Gerstenhaber \cite{Gerstenhaber63} on Hochschild cohomology and
deformations of algebras, under the name ``pre-Lie algebras'' which is now the
standard terminology. The term \textsl{chronological algebras} has also been
sometimes used, e.g. in \cite{AG81}. Antisymmetrising the pre-Lie product gives rise to a Lie algebra, the universal enveloping algebra $\mathcal U(P)$ of which is isomorphic (as a Hopf algebra) to the symmetric algebra $S(P)$ endowed with the deshuffle coproduct
\[\Delta_{\sshu}(x_1\cdots x_n)=\sum_{I\sqcup J=\{1,\ldots, n\}}x_I\otimes x_J\]
(with the notation $x_J:=x_{j_1}\cdots x_{j_k}$ for $J=\{j_1,\ldots, j_k\}$) and  the Grossman-Larson product
\[X\star Y:=\sum_{(X)}X_1(X_2\rhd Y),\]
where the Guin-Oudom product $\rhd$  \cite[Proposition 2.7]{GO08} is the unique extension of the pre-Lie product to $S(P)$ such that
\begin{itemize}
\item $1\rhd X=X$,
\item $X\rhd YZ= \sum_{(X)} (X_{1} \rhd Y) (X_{2} \rhd Z)$,
\item $(xY)\rhd Z=x\rhd(Y\rhd Z)-(x\rhd Y)\rhd Z$
\end{itemize}
for any $X,Y,Z\in S(L)$ and $x\in L$. A key property of the Grossman-Larson product is given by
\[X\rhd(Y\rhd Z)=(X\star Y)\rhd Z\]
for any $X,Y,Z\in S(P)$. For a short survey on pre-lie algebras, see \cite{Bu06,M11}.\\

F. Chapoton and M. Livernet provided the
explicit description of the free pre-Lie algebra $\mop{PL}(A)$ in terms of $A$-decorated rooted
trees endowed with grafting \cite{CL01}. Recall that a rooted tree is a connected oriented graph with a finite number of vertices,
one among them being distinguished as the root, such that any vertex
admits exactly one outgoing edge, except the root which has no outgoing
edges. Here is the list of rooted trees up to five vertices, where the edges are tacitly oriented from top to bottom:
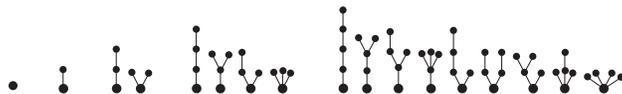
\begin{equation*}
\racine \hskip 5mm \arbrea \hskip 5mm  \arbreba \arbrebb \hskip 5mm  \arbreca \arbrecb \arbrecc
\arbrecd \hskip 5mm  \arbreda \arbredb \arbredc \arbredd \arbrede \arbredf \arbredz \arbredg
\arbredh
\end{equation*}
A rooted forest is a finite collection of rooted trees, possibly with repetitions. For any set $A$, an $A$-decorated rooted forest $f$ is a rooted forest together with a map $d:\mathcal V(f)\to A$, where $\mathcal V(f)$ is the set of vertices of $f$. The pre-Lie product $s\to t$ of two rooted trees is obtained by grafting the root of $s$ on a vertex of $t$, summing up over all choices of vertices:
\begin{equation}\label{greffe}
s\to t=\sum_{v\in\mathcal V(t)}s\to_v t.
\end{equation}
For example,
\[\racine\to\arbrea=\arbrebb+\arbreba.\]
The symmetric algebra $S\big(\mop{PL}(A)\big)$ is the linear span of $A$-decorated rooted forests. Denoting also by $\to$ the Guin-Oudom extension of the grafting product, the corresponding Grossman-Larson product \cite{GL89} is given by the \textsl{graft-or-fall} formula
\[F\star G=B_-\big(F\to B_+^c(G)\big).\]
Here $B_+^c(F):=F\to\racine_c$ for some $c\in A$, and $B_-$ is the root removal, left inverse to $B_+^c$.\\

The Butcher--Connes--Kreimer Hopf algebra (\cite{CK98, F02}, see also \cite{B63}) $\mathcal
  H^A_{\smop{BCK}}=\bigoplus_{n\geq 0} \left(\mathcal
  H^A_{\smop{BCK}}\right)_n$ is a commutative Hopf algebra of $A$-decorated rooted forests over
$\mathbb K$, graded by the number of vertices, obtained as the graded dual of the cocommutative Grossman-Larson Hopf algebra described above. Normalising the dual forest basis by the symmetry factor (see e.g. \cite{Br00, MS08}), namely
\begin{equation}\label{sym-tree}
\langle u,v\rangle :=\sigma(u)\delta_u^v
\end{equation}
where $\sigma(u)$ is the cardinal of the automorphism group of the forest $u$, the coproduct on a rooted forest $u$ is described as follows: the set $\mathcal V(u)$ of vertices of a forest $u$ is endowed with a partial order defined by $x \le y$ if and only if there is a path from $y$ to a root passing through $x$. Any subset $W$ of $\mathcal V(u)$ defines a subforest $u\restr W$ of $u$ in an obvious manner, i.e. by keeping the edges of $u$ which link two elements of $W$. The coproduct is then defined by:
\begin{equation}
\label{coprod}
	\Delta_{\smop{BCK}}(u)= \sum_{V \amalg W=\mathcal V(u) \atop W<V}u\restr V\otimes u\restr W.
\end{equation}
Here the notation $W<V$ means that $y\not >x$ for any vertex $x$ in $V$ and any
vertex $y$ in $W$. Such a couple $(V,W)$
is also called an {\sl admissible cut\/}, with crown (or pruning) $u\restr V$ and
trunk $u\restr W$. We have for example:
 \allowdisplaybreaks{
\begin{eqnarray*}
\Delta_{\smop{BCK}}\big(\arbrea\big) &=&
                        \arbrea \otimes \un
                            + \un \otimes \arbrea
                               + \racine \otimes \racine \\
 \Delta_{\smop{BCK}}\big(\! \arbrebb \big) &=&
 					\arbrebb \otimes \un
                               + \un \otimes \arbrebb +
                                  2\racine \otimes\arbrea
                                   + \racine\racine\otimes \racine \, .
\end{eqnarray*}}

\noindent The iterated restricted coproduct is best written in terms of ordered weak partitions of $\mathcal V(u)$ into $n$ blocks:
\begin{equation}
\label{iter-coprod}
\wt\Delta_{\smop{BCK}}^{n-1}(u)=\sum_{V_1\amalg\cdots\amalg V_n=\mathcal V(u) \atop V_n<\cdots <V_1}u\restr{V_1}\otimes\cdots\otimes u\restr{V_n}\,.
\end{equation}
Coassociativity of the coproduct follows immediately. Note however that the relation $<$ on subsets of vertices is not transitive: the notation $V_n<\cdots <V_1$ is to be understood as $V_i<V_j$ for any $i,j\in\{1,\ldots ,n\}$ with $i>j$.\\

Recall that the symmetry factor of a rooted forest is the cardinal of its automorphism group. For later use, the symmetry factor of a forest $u=t_1^{\ell_1}\cdots t_n^{\ell_n}$, where $t_1,\ldots, t_n$ are pairwise distinct trees, can be expressed in terms of the symmetry factor of its tree components as
\begin{equation}\label{sym-forest}
\sigma(u)=\ell_1!\cdots\ell_n!\sigma(t_1)^{\ell_1}\cdots\sigma(t_n)^{\ell_n}.
\end{equation}
\begin{rmk}\label{int-ext-forests}\rm
The symmetry factor $\sigma (u)$ is obtained by multiplying the \textsl{external symmetry factor} $\sigma^{\smop{ext}} (u):=\ell_1!\cdots \ell_n!$ by the \textsl{internal symmetry factor} $\sigma^{\smop{int}} (u):=\sigma(t_1)^{\ell_1}\cdots\sigma(t_n)^{\ell_n}$. The internal symmetry factor is the cardinal of the internal automorphism group, i.e. the normal subgroup of automorphisms which leaves each connected component unchanged. The external symmetry factor is the cardinal of the external automorphism group, which is a subgroup of the permutation group of the connected components.
\end{rmk}
\noindent Let us finally recall that the coproduct can be equivalently written as
\begin{equation}\label{bck-adm}
\Delta_{\smop{BCK}}(u)=\sum_{c\in\smop{Adm}(u)}P^c(u)\otimes R^c(u).
\end{equation}
\ignore{
\begin{defn}
Let $v$ and $w$ be two rooted forests, with $v=t_1\cdots t_r$, where the $t_j$'s are the connected components of $v$ (not necessarily pairwise disjoint). A grafting of $v$ on $w$ is a rooted forest $u$ obtained by choosing a collection $(x_j)_{j=1,\ldots, r}$ of vertices of $w$ (possibly with repetitions), and by grafting each tree $t_j$ on the corresponding vertex $x_j$.
\end{defn}
\begin{prop}\label{sym-cut-graft}
Let $u,v,w$ be three decorated rooted forests. Let $C(u,v,w)$ be the number of admissible cuts of $u$ such that $P^c(u)=v$ and $R^c(u)=w$. Let $G(u,v,w)$ the number of graftings of $v$ on $w$ isomorphic to $u$. Then
\[C(u,v,w)=\frac{\sigma(u)}{\sigma(v)\sigma(w)}G(u,v,w).\]
\end{prop}
\begin{proof}
Let $\mathcal C(u,v,w)$ be the set of admissible cuts of $u$ such that $P^c(u)=v$ and $R^c(u)=w$, and let $\mathcal G(u,v,w)$ be the set of graftings of $v$ on $w$ isomorphic to $u$. We denote by $U, V, W$ the automorphism groups of $u$, $v$ and $w$ respectively. For any subset $E$ of vertices of $w$, we denote by $W_E$ the subgroup of $W$ which stabilises any element of $E$.\\

For any $c\in\mathcal C(u,v,w)$, we denote by $E'_c$ the set of vertices below the cut (i.e. the maxima of $R^c(u)$ admitting at least one vertex of $P^c(u)$ above them). For any $x\in E'_c$, we denote by $F'_x$ the set of vertices $y$ of $P^c(u)$ such that there is an edge from $x$ to $y$. Let $W'_c$ be the subgroup of $\mop{Aut} R^c(u)$ which stabilises any element of $E'_c$, and let $V'_c$ be the subgroup of $\mop{Aut} P^c(u)$ which stabilises the subsets $F'_x$ for any $x\in E'_c$. The group $U$ acts transitively on $\mathcal C(u,v,w)$, and the stabiliser of $c$ is $W'_c\times V'_c$.\\

On the other hand, for any grafting $b\in\mathcal G(u,v,w)$, let $E_b$ be the set of vertices of $w$ on which the grafting occurs, and, for any $x\in E_b$, let $F_x$ be the set of roots of $v$ which are grafted on $x$.  Let $W_b$ be the subgroup of $W$ which stabilises any element of $E_b$, and let $V_b$ be the subgroup of $V$ which stabilises the subsets $F_x$ for any $x\in E_c$. The group $W\times V$ acts transitively on $\mathcal G(u,v,w)$, and the stabiliser of $b$ is $W_b\times V_b$. By the orbit-stabiliser theorem we have
\begin{equation}\label{os}
C(u,v,w)=\frac{|U|}{|W'_c|. |V'_c|}, \hskip 8mm G(u,v,w)=\frac{|V|.|W|}{|W_b|.|V_b|}.
\end{equation}
For any grafting $b$, if we denote by $c$ the admissible cut which "undoes" the grafting, it is easily seen that the groups $W'_c$ and $W_b$ are isomorphic, and similarly $V'_c$ and $V_b$. From \eqref{os} we therefore get
\[\frac{C(u,v,w)}{G(u,v,w)}=\frac{|U|}{|V|.|W|},\]
which ends up the proof.
\end{proof}
}
\subsection{Novikov algebras and the Linares--Otto--Tempelmayr Hopf algebra of multi-indices}\label{par:nov}
A Novikov algebra  \cite{DL02} is a vector space $N$ over a base field $\mathbb K$, together with a bilinear product $\rhd:N\times N\to N$ such that, for any $x,y,z\in N$, the following identities hold:
\begin{align}
x\rhd(y\rhd z)-(x\rhd y)\rhd z &= y\rhd(x\rhd z)-(y\rhd x)\rhd z,\\
(x\rhd y)\rhd z&=(x\rhd z)\rhd y.
\end{align}
The first is left pre-Lie identity, the second is right-NAP identity\footnote{NAP for Non-Associative Permutative \cite{L06}. Novikov algebras are right pre-Lie and left NAP in the definition given in \cite{DL02} (right-symmetric and left-commutative in the terminology employed therein).}. Novikov algebras seem to appear for the first time in the article \cite{GD79} by I. M. Gelfand and I. J. Dorfman in the study of Hamiltonian operators in the formal calculus of variations (see Equation (6.3) therein, where the left pre-Lie identity appears in a disguised form). They have been rediscovered by A. A. Balinskii and S. P. Novikov in \cite{BN85} (in the right pre-Lie and left NAP form, see Equation (4) therein). The terminology was proposed by J. M. Osborn in \cite{O92}. An important example is given by a commutative associative algebra $\mathcal A$ endowed with a derivation $D:\mathcal A\to \mathcal A$. From Leibniz' rule $D(xy)=D(x)y+xD(y)$, it is easily seen that $(\mathcal A,\rhd)$ is Novikov with $x\rhd y:=xD(y)$.\\

With our conventions, the free Novikov algebra $\mop N(A)$ generated by a set $A$ is described as follows \cite[Definition 7.7]{DL02}: let $\overline{\mop N}(A)$ be the commutative algebra of polynomials with variables $x_j^a,\, (a,j)\in A\times \{-1,0,1,2,\ldots\}$, let $\partial$ be the unique derivation of $\overline{\mop N}(A)$ such that $\partial x_j^a=x_{j+1}^a$. The integer $j\ge -1$ is the \textsl{weight} of the variable $x_j^a$.\\

\noindent A basis of $\overline{\mop N}(A)$ is given by the monomials
\[\mathbf x^{\mathbf k}:=\prod_{j\ge -1,\, a\in A}(x_j^a)^{k_j^a},\]
where the exponents $k_j^a$ are non-negative integers, equal to zero except a finite number of them. The weight induces a unique $\mathbb Z$-grading of the algebra $\overline{\mop N}(A)$, for which the derivation $\partial$ is homogeneous of degree one. By Leibniz' rule, the expression of the derivation $\partial$ is given by
\begin{equation}\label{expression-d}
\partial\mathbf x^{\mathbf k}=\sum_{j\ge -1,\,a\in A}k_j^a\mathbf x^{\mathbf k-\mathbf e_j^a+\mathbf e_{j+1}^a},
\end{equation}
where $\mathbf e_j^a$ is the multi-index in which all coordinates are equal to zero except the one in position $(j,a)$ which is equal to one. In other words, $\mathbf x^{\mathbf e_j^a}=x_j^a$.
The bilinear product $P\rhd Q:=P.\partial Q$ endows $\overline{\mop N}(A)$ with a Novikov structure, and the free Novikov algebra $\mop N(A)$ turns out to be the homogeneous component of $\overline{\mop N}(A)$ of weight $-1$. This is Theorem 7.8 in \cite{DL02}, recently generalised to multi-Novikov algebras \cite[Theorem 3.3]{BD2023}. The embedding of $A$ into $\mop N(A)$ is given by $a\mapsto x_{-1}^a$.\\

\noindent We shall also use the unique derivation $\overline \partial:\overline N(A)\to\overline N(A)$ defined on the variables by
\begin{equation}\label{dbar}
\overline \partial x_j^a=x_{j-1}^a \hbox{ if } j\ge 0,\hskip 8mm \overline \partial x_{-1}^a=0.
\end{equation}
Its expression on the basis of monomials is given by
\begin{equation}\label{expression-dbar}
\overline\partial\mathbf x^{\mathbf k}=\sum_{j\ge 0,\,a\in A}k_j^a\mathbf x^{\mathbf k-\mathbf e_j^a+\mathbf e_{j-1}^a}.
\end{equation}
The canonical surjective map $\Phi$ from $\mop{PL}(A)$ to $\mop{N}(A)$ is the unique pre-Lie morphism extending the embedding of the generators, and can be understood as the \textsl{fertility map}: for any $A$-decorated rooted tree $t$, we have
\begin{equation}\label{F-map}
\Phi(t)=\prod_{v\in\mathcal V(t)}x^{d(v)}_{f(v)-1}\,,
\end{equation}
where $d(v)\in A$ is the decoration of vertex $v$, and where $f(v)$ is its fertility, i.e. its number of incoming edges. The map $\Phi$ is obviously surjective. From \eqref{F-map} and \eqref{greffe}, it is easily seen to be a pre-Lie morphism, namely
\begin{equation}\label{Fplm}
\Phi(s\to t)=\Phi(s).\partial \Phi(t)
\end{equation}
for any rooted trees $s,t$, due to the fact that grafting any tree at vertex $v$ increases its fertility by $1$, leaving the fertility of the other vertices unchanged\footnote{The map $\Phi$ appears in the recent literature, e.g. \cite[Paragraph 2.5]{BL2023}, with extra symmetry factors, due to a different normalisation of the derivation $\partial$. See Remark \ref{alt-sym} below.}.

\begin{rmk}\rm
Recall that $C^\infty(\mathbb R^d,\mathbb R^d)$, identified with the set of vector fields on $\mathbb R^d$, is a pre-Lie algebra with pre-Lie product
\[\left(\sum_{i=1}^d f_i\partial_i\right)\rhd \left(\sum_{j=1}^d f_j\partial_j\right)=\sum_{i=1}^df_i\left(\sum_{j=1}^d\partial_if_j\right)\partial_j.\]
When $d=1$, this product specialises to the Novikov product
\[f\rhd g=f.\partial g\]
where $\partial g$ stands for the derivative of $g$. By freeness universal property, for any family $f:=(f_a)_{a\in A}$ of vector fields on $\mathbb R$, there is a unique Novikov algebra morphism $\mathcal F_f:\mop N(A)\to C^\infty(\mathbb R,\mathbb R)$ such that $\mathcal F_f(\mathbf x_{-1}^a)=f_a$.
This is the starting point of multi-index B-series introduced in \cite{BEH2024}:
\[B_f(\alpha)=\sum_{\mathbf k,\, \smop{wt}(\mathbf x^{\mathbf k})=-1}\frac{\alpha(\mathbf x^{\mathbf k})}{\sigma(\mathbf x^{\mathbf k})}\mathcal F_f(\mathbf x^{\mathbf k}),\]
where $\alpha$ is any linear map from $\mop N(A)$ into $\mathbb R$. Our convention for the symmetry factor $\sigma(\mathbf x^{\mathbf k})$ is given in Paragraph \ref{par:pair} below.
\end{rmk}

As to any pre-Lie algebra, the Guin-Oudom procedure applies to both $\mop{N}(A)$ and $\mop{PL}(A)$. In particular, the multiplicative extension
\[\Phi:S\big(\mop{PL}(A)\big)\longrightarrow S\big(\mop{N}(A)\big)\]
of the fertility map to symmetric algebras is a Hopf algebra morphism from $\Big(S\big(\mop {PL}(A)\big),\star, \Delta_{\sshu},\mathbf 1,\varepsilon\Big)$ onto $\Big(S\big(\mop N(A)\big),\star, \Delta_{\sshu},\mathbf 1,\varepsilon\Big)$, where both Grossman-Larson products are denoted by $\star$, and where $\Delta_{\sshu}$ is the usual deshuffle coproduct. The Hopf algebra of multi-indices is defined by
\[\mathcal H_{\smop{LOT}}^A=\Big(S\big(\mop N(A)\big),\star, \Delta_{\sshu},\mathbf 1,\varepsilon\Big)^\circ,\]
where $(-)^\circ$ stands for the graded dual. A basis of $S\big(\mop N(A)\big)$ is given by \textsl{monomials of monomials} $\mathbb M=M_1^{m_1}\odot\cdots\odot M_k^{m_k}$ where the $M_j$'s are distinct monomials in $\mop N(A)$ (we use the notation $\odot$ for the commutative ``external'' product in the symmetric algebra, not to be confused with the commutative ``internal'' product of $\overline N(A)$). The explicit formula for the coproduct depends on the choice of a pairing. A proposal will be given in the next section.
\section{Pairing and embedding}\label{sect:pair}
\subsection{The pairing}\label{par:pair}
\noindent The symmetry factor of any monomial
\[M=\mathbf x^{\mathbf k}=\prod_{a\in A,\,j\ge -1}(x_j^a)^{k_j^a},\]
where $\mathbf k$ stands for the multi-index $(k_j^a)_{a\in A,\, j\ge -1}$, is given by
\[\sigma(M)=\mathbf k!=\prod_{a\in A,\,j\ge -1}k_j^a!\,.\]
The degree of the monomial $\mathbf x^{\mathbf k}$ is
\[|\mathbf k|:=\sum_{a\in A,\,j\ge -1}k_j^a.\]
Its weight is defined by
\[\mop{wt}(\mathbf x^{\mathbf k}):=\sum_{a\in A,\,j\ge -1}jk_j^a,\]
which we also will denote by $\mop{wt}(\mathbf k)$. Now let us define a symmetric nondegenerate pairing on $S\big(\overline{\mop{N}}(A)\big)$ as follows: let $\mathbb M=M_1^{\odot \ell_1}\odot\cdots\odot M_n^{\odot \ell_n}$ be a monomial of monomials, and let $\mathbb M'$ be another one, and set
\[\langle \mathbb M,\,\mathbb M'\rangle:=\sigma(\mathbb M)\delta_{\mathbb M}^{\mathbb M'},\]
where the symmetry factor is given by a formula similar to \eqref{sym-forest}:
\begin{equation}
\sigma (\mathbb M):=\ell_1!\cdots \ell_n! \sigma(M_1)^{\ell_1}\cdots\sigma(M_n)^{\ell_n}.
\end{equation}
Under this pairing, the dual of the deshuffle coproduct is the commutative product $\odot$. It remains to compute the dual coproduct of the Grossman-Larson product.
\begin{rmk}\label{int-ext}\rm
The symmetry factor $\sigma (\mathbb M)$ is obtained by multiplying the \textsl{external symmetry factor} $\sigma^{\smop{ext}} (\mathbb M):=\ell_1!\cdots \ell_n!$ by the \textsl{internal symmetry factor} $\sigma^{\smop{int}} (\mathbb M):=\sigma(M_1)^{\ell_1}\cdots\sigma(M_n)^{\ell_n}$.
\end{rmk}
\begin{prop}\label{ddbar-duality}
Considering the restriction of the above pairing to the polynomial algebra $\overline{\mop N}(A)$, the derivation $\overline\partial$ is the transpose of the derivation $\partial$.
\end{prop}
\begin{proof}
By direct computation:
\[
\begin{split}
\langle \partial\mathbf x^{\mathbf k},\mathbf x^\ell \rangle& = \sum_{j\geq -1,\, a\in A} k_j^a \langle \mathbf x^{\mathbf k+\mathbf e^a_{j+1}-\mathbf e^a_j},\mathbf x^\ell \rangle 
\\ &=\sum_{j\geq -1, \, a\in A} k_j^a \,\ell! \, \ind_{(\mathbf k+\mathbf e^a_{j+1}-\mathbf e^a_j=\ell)}
\\ & = \sum_{j\geq -1,\, a\in A} (\ell+\mathbf e^a_j)! \, \ind_{(\mathbf k+\mathbf e^a_{j+1}=\ell+\mathbf e^a_j)}.
\end{split}
\]
Now
\[
\begin{split}
\langle \mathbf x^{\mathbf k},\bar\partial\mathbf x^\ell \rangle& = \sum_{j\geq 0,\, a\in A} \ell^a_j \langle \mathbf x^{\mathbf k},\mathbf x^{\ell+\mathbf e^a_{j-1}-\mathbf e^a_j} \rangle
\\ &=\sum_{j\geq 0,\, a\in A} \ell^a_j (\ell+\mathbf e^a_{j-1}-\mathbf e^a_j)! \, \ind_{(\mathbf k=\ell+\mathbf e^a_{j-1}-\mathbf e^a_j)}
\\ & = \sum_{j\geq 0, \, a\in A} (\ell+\mathbf e^a_{j-1})! \, \ind_{(\mathbf k+\mathbf e^a_{j}=\ell+\mathbf e^a_{j-1})},
\end{split}
\]
hence both expressions coincide.
\end{proof}
For any multi-index $\mathbf k$, let us denote by $\shift{\mathbf k}$ the corresponding left-shifted multi-index, with coordinates ${\shift k}_j^a:=k_{j+1}^a$. In particular, we have
\[\shift{\mathbf e_j^a}=\mathbf e_{j-1}^a.\]
\begin{prop}\label{dbar-r}
For any integer $r\ge 0$ and for any multi-index $\mathbf k$ we have
\begin{equation}\label{ckl}
\bar\partial^r \mathbf x^{\mathbf k} = \sum_{|\ell|=r} C_{\mathbf k,\ell}\mathbf x^{\mathbf k-\ell+\shift\ell},
\end{equation}
where the coefficients $C_{\mathbf k,\ell}$ are given by $C_{\mathbf k,0}=1$ and the recursive formula below:
\begin{equation}\label{ckl-rec}
C_{\mathbf k,\ell}=\sum_{j\ge -1,\, a\in A,\, \ell_j\ge 1}C_{\mathbf k,\ell-\mathbf e^a_j}(k^a_j-\ell_j^a+1+\ell^a_{j+1}).
\end{equation}
\end{prop}
\begin{proof}
The definition of $\overline\partial$ immediately yields
\begin{equation}
C_{\mathbf k,\mathbf e_j^a}=k_j^a.
\end{equation}
We therefore can compute:
\begin{eqnarray*}
\overline\partial^{r+1}\mathbf x^{\mathbf k}&=&\sum_{|\ell'|=r}C_{\mathbf k,\ell'}\mathbf x^{\mathbf k+\shift{\ell'}-\ell'}\\
&=&\sum_{|\ell'|=r}\sum_{j\ge 0,\, a\in A}C_{\mathbf k,\ell'}(k_j^a+\ell_{j+1}^a-{\ell'_j}^a)\mathbf x^{\mathbf k+\shift{(\ell'+\mathbf e_j^a)}-(\ell'+\mathbf e_j^a)}\\
&=&\sum_{|\ell|=r+1}\sum_{(\ell',j,a),\; \ell'+\mathbf e^a_j=\ell}C_{\mathbf k,\ell-\mathbf e_j^a}(k_j^a+\ell_{j+1}^a-\ell_j^a+1) \mathbf x^{\mathbf k+\shift\ell -\ell},
\end{eqnarray*}
which proves the claim.
\end{proof}
\ignore{
\noindent For later use, we have the following explicit expression: let $r$ be a non-negative integer and let $\mathbf x^{\mathbf k}$ be a monomial. Then
\begin{equation}\label{dbar-r}
\overline\partial^r(\mathbf x^{\mathbf k})=\sum_{|\ell|=r}\frac{\mathbf k!}{(\mathbf k-\mathbf \ell)!}\mathbf x^{\mathbf k-\mathbf \ell}.
\end{equation}
}
\begin{prop}\label{prop:sym-forest}
Let $(M_1,\ldots, M_r)$ be an $r$-tuple of monomials in $N(A)$, and let $\mathbb M=M_1\odot\cdots\odot M_r$ be the monomial of monomials obtained by multiplying the $M_j$'s together. Let $F$ be a rooted decorated forest such that $\Phi(F)=\mathbb M$, and let $\mathcal A$ be the set of $r$-tuples of decorated rooted trees given by
\[\mathcal A:=\{(t_1,\ldots, t_r), \, t_1\cdots t_r=F \hbox{ and }\Phi(t_j)=M_j \hbox{ for any }j=1,\ldots, r\}.\]Then we have
\[\vert \mathcal A\vert = \frac{\sigma^{\smop{ext}}(\mathbb M)}{\sigma^{\smop{ext}}(F)}.\]
\end{prop}
\begin{proof}
The external automorphism group of $\mathbb M$ acts transitively on $\mathcal A$. The stabiliser of the $r$-tuple $(t_1,\ldots, t_r)$ is the external automorphism group of the forest $F$. One concludes by the orbit-stabiliser theorem.
\end{proof}
\begin{rmk}\label{alt-sym}\rm
An alternative choice of symmetry factors is frequently used in the literature: in \cite{BH24} (see also \cite[Lemma 6.1]{LOT23} and \cite[Section 2]{BEH2024}), the authors use (translated in our notations and context)
\[\wt\sigma(\mathbf x^{\mathbf k})=\prod_{a\in A,\, j\ge -1}\big((j+1)!\big)^{k_j^a}.\]
Considering that a variable $x_j^a$ corresponds to a vertex decorated by $a$ with $j+1$ edges above it, one may consider that it has an intrinsic symmetry factor $(j+1)!$. For a multi-index $\mathbf k=(k_j^a)_{a\in A,\, j\ge -1}$, the symmetry factor above can therefore be considered as the internal part of a symmetry factor
\[\wh\sigma(\mathbf x^{\mathbf k}):=\sigma(\mathbf x^{\mathbf k})\wt\sigma(\mathbf x^{\mathbf k})=\mathbf k!\prod_{a\in A,\, j\ge -1}\big((j+1)!\big)^{k_j^a},\]
our $\mathbf k!$ being the external part. Note also that our definition of the derivation $\partial$ on $\overline{\mop{N}}(A)$ differs from \cite[Paragraph 3.2]{LOT23}, where the authors choose
\[\partial x_j^a=(j+1)x_{j+1}^a.\]

\end{rmk}
\subsection{The embedding}\label{par:jmath}
Define the embedding $\jmath:\mathcal H_{\smop{LOT}}^A\inj 6\mathcal H_{\smop{BCK}}^A$ by transposing the multiplicatively extended fertility map $\Phi$. This is an injective Hopf algebra morphism. From $\langle \Phi(t),\mathbf x^{\mathbf k}\rangle=\langle t,\jmath(\mathbf x^{\mathbf k})\rangle$ for any rooted tree $t$ and any monomial $\mathbf x^{\mathbf k}$ of weight $-1$, and considering the definitions of both pairings, we easily get
\begin{equation}\label{jmath}
\jmath(\mathbf x^{\mathbf k})=\sum_{t,\, \Phi(t)=\mathbf x^{\mathbf k}}\frac{\sigma(\mathbf x^{\mathbf k})}{\sigma(t)}t.
\end{equation}
For example,
\[\jmath(x_{-1}x_0)=\arbrea,\hskip 8mm \jmath(x_{-1}^2x_1)=\arbrebb,\hskip 8mm \jmath(x_{-1}^2x_0x_1)=2\arbrecc+\arbrecb, \hskip 8mm\jmath(x_{-1}^3x_0x_2)=\arbredd+3\arbredg .\]
\begin{lem}\label{prop-jdbar-r}
\begin{equation}\label{jdbar-r}
\jmath(\overline\partial^r\mathbf x^{\mathbf k})=\sum_{|\ell|=r}\,\sum_{t,\, \Phi(t)=\mathbf x^{\mathbf k+\shift\ell-\ell}}\frac{D_{\mathbf k,\ell}}{\sigma(t)}t,
\end{equation}
with $D_{\mathbf k,\ell}=C_{\mathbf k,\ell}(\mathbf k+\shift\ell-\ell)!$. The coefficients $D_{\mathbf k,\ell}$ are given by $D_{\mathbf k,0}=\mathbf k!$ and the recursive formula
\begin{equation}\label{jdbar-r-rec}
D_{\mathbf k,\ell}=\sum_{j\ge -1,\, a\in A, \ell_j^a\ge 1}D_{\mathbf k,\ell-\mathbf e_j^a}(k_{j-1}^a-\ell_{j-1}^a+\ell_j^a).
\end{equation}
\end{lem}
\begin{proof}
Equation \eqref{jdbar-r} is a direct consequence of \eqref{ckl} and \eqref{jmath}. From \eqref{ckl-rec} and the expression of $D_{\mathbf k,\ell}$ in terms of $C_{\mathbf k,\ell}$, we can compute:
\begin{eqnarray*}
D_{\mathbf k,\ell}&=&\sum_{j\ge-1,\, a\in A,\, \ell_j^a\ge 1}C_{\mathbf k,\ell-\mathbf e_j^a}(k_j^a-\ell_j^a+1+\ell_{j+1}^a)(\mathbf k+\shift\ell-\ell)!\\
&=&\sum_{j\ge-1,\, a\in A,\, \ell_j^a\ge 1}D_{\mathbf k,\ell-\mathbf e_j^a}\frac{(\mathbf k+\shift\ell-\ell)!}{(\mathbf k+\shift\ell-\mathbf e_{j-1}^a-\ell+\mathbf e_j^a)!}(k_j^a-\ell_j^a+1+\ell_{j+1}^a)\\
&=&\sum_{j\ge -1,\, a\in A, \ell_j^a\ge 1}D_{\mathbf k,\ell-\mathbf e_j^a}(k_{j-1}^a-\ell_{j-1}^a+\ell_j^a).
\end{eqnarray*}
\end{proof}
\subsection{Cocycles and mock-cocycles}\label{par:bplus}
Recall the cocycle operator $B_+^a:S\big(\mop{PL}(A)\big)\to\mop{PL}(A)$ for any $a\in A$, which grafts all trees of a given forest on a common root decorated by $a$. Any $A$-decorated rooted tree $t$ is uniquely given by $t=B_+^a (F)$, where $F$ is the forest obtained from $t$ by removing the root, and where $a$ is the decoration of the root of $t$. Its transpose $B_-^a:\mathcal H_{\smop{BCK}}^A\to \mathcal H_{\smop{BCK}}^A$ is given by $B_-^a(t)= F$ if $t=B_+^a(F)$, and by $B_-^a(t)= 0$ if the decoration of the root of $t$ is different from $a$.\\

\noindent Now consider, for any $a\in A$, the operator $L^a:S\big(\mop{N}(A)\big)\to\mop N(A)$ defined by
\[L^a(\mathbf x^{\mathbf k^1}\odot\cdots\odot \mathbf x^{\mathbf k^r}):=\mathbf x^{\mathbf k^1+\cdots+\mathbf k^r}x_{r-1}^a.\]
We obviously have
\begin{equation}\label{phib}
\Phi\circ B_+^a=L^a\circ \Phi.
\end{equation}
Its transpose $\overline L^a:\mathcal H_{\smop{LOT}}^A\to \mathcal H_{\smop{LOT}}^A$ is given by
\[\overline L^a(\mathbf x^{\mathbf k})=\sum_{\mathbb M,\, L^a(\mathbb M)=\mathbf x^{\mathbf k}}\frac{\sigma(\mathbf x^{\mathbf k})}{\sigma(\mathbb M)}\mathbb M.\]
From \eqref{phib} we immediately get
\begin{equation}\label{phibtr}
\jmath\circ \overline L^a=B_-^a\circ\jmath.
\end{equation}
\section{Free edges, graftings and cuts}\label{sect:bc}
\subsection{Rooted trees and forests with free edges}
The free Novikov $N(A)$ is included inside a wider Novikov algebra $\overline N(A)$. Similarly, the free pre-Lie algebra $\mop{PL}(A)$ is included in a wider pre-Lie algebra $\overline{\mop{PL}}(A)$ defined as the linear span of rooted trees with free edges, i.e. edges without upper vertex. It can be seen as the free pre-Lie algebra generated by $A\times \mathbb N_0$, where the second component of the decoration of a given vertex indicates the number of free edges attached to it.\\

The linear span of rooted forests with free edges will be denoted by $\overline {\mathcal H}_{\smop{BCK}}^A$.\ignore{The grafting of two such trees $s$ and $t$ is still defined by \eqref{greffe}.} The weight $\mop{wt}(t)$ of $t$ is given by the total number of edges minus the number of vertices. The pairing is naturally extended to  $\overline {\mathcal H}_{\smop{BCK}}^A$, where the symmetry factor of a forest with free edges is understood as the symmetry factor of the corresponding $A\times\mathbb N_0$-decorated forest. We obviously have $\mop{wt}(t)=-1$ for an ordinary rooted tree without free edges.
\[\arbrefr\]
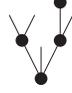
\centerline{\small A rooted tree of weight $3$, with four free edges.}\\

\noindent Formula \eqref{F-map} extends the fertility map $\Phi$ to a map $\overline \Phi:\overline{\mop{PL}}(A)\to\overline N(A)$ which respects the weight. Similarly, Formula \eqref{jmath} extends the embedding $\jmath$ to a map
\[\overline\jmath:\overline {\mathcal H}_{\smop{LOT}}^A\longrightarrow \overline {\mathcal H}_{\smop{BCK}}^A,\]
where $\overline {\mathcal H}_{\smop{LOT}}^A$ stands for the commutative algebra $S\big(\overline{N}(A)\big)^\circ$. We remark that $\overline\jmath$ is not an embedding, as it sends any monomial $\mathbf x^{\mathbf k}$ of weight $\le -2$ to zero.\\

\noindent Let $\delta: \overline{\mop{PL}}(A)\to \overline{\mop{PL}}(A)$ be the map given by
\[\delta(t)=\sum_{v\in\mathcal V(t)}\delta_v(t),\]
where $\delta_v$ adds a free edge to vertex $v$. For example,
\[\delta \Big(\arbrefrun\Big)= \arbrefrpetit + \arbrefrpetitb + \arbrefrpetitc + \arbrefrpetitd.\]
Let us use the symbol $\overline\delta$ for the map given by
\[\overline\delta(t)=\sum_{v\in\mathcal V(t)}\overline\delta_v(t),\]
where $\overline\delta_v$ removes a free edge at vertex $v$, and returns zero if $v$ has no free edge. For example,
\[\overline\delta \Big(\arbrefrpetit\Big)=\arbrefrun + \arbrefrdeux + \arbrefrtrois.\]
\begin{prop}
The map $\overline\delta$ is the transpose of $\delta$, and we have
\begin{equation}
\overline\Phi\circ\delta=\partial\circ\overline \Phi.
\end{equation}
\end{prop}
\begin{proof}
Choosing a vertex $v$ (resp. $w$) in a tree $t$ (resp. $u$), and
considering both sets
\[P:=\{v\in \mathcal V(t),\, \delta_v(t)=u\},\hskip 12mm Q:=\{w\in\mathcal V(u),\, \overline\delta_w(u)=t\},\]
the orbit-stabiliser theorem states
\[P\simeq \mop{Aut} (t)/\mop{Aut}_v (t),\hskip 12mm Q\simeq \mop{Aut} (u)/\mop{Aut}_w (u),\]
where $v$ (resp. $w$) is any choice of element in $P$ (resp. Q), and where $\mop{Aut}_v (t)$ is the stabiliser of $v$ in  $\mop{Aut} (t)$ (and similarly for $Q$). If two trees $t$ and $u$ are such that $u=\delta_v(t)$ for some $v\in\mathcal V(t)$, then both sets of vertices $\mathcal V(t)$ and $\mathcal V(u)$ can be naturally identified. In that case, from the obvious isomorphism $\mop{Aut}_v (t)\simeq \mop{Aut}_v (u)$ we get
\[\sigma(u)|P|=\sigma(t)|Q|.\]
We therefore have
\begin{eqnarray*}
\langle \delta(t),u\rangle &=&|P|\sigma(u)\\
&=&|Q|\sigma(t)\\
&=&\langle t,\overline\delta(u)\rangle,
\end{eqnarray*}
which proves the first assertion. The second assertion comes from the fact that adding a free edge to a vertex corresponds to shifting the variable associated to it. 
\end{proof}
\begin{cor}\label{deltabar-j}
\begin{equation}
\overline\delta\circ\overline\jmath=\overline\jmath\circ\overline\partial.
\end{equation}
\end{cor}
\begin{prop}\label{deltabar-max}
Let $t$ be a rooted $A$-decorated tree with $r$ free edges. Let $r_v$ be the number of free edges at a given vertex $v$ of $t$. Then we have
\[\overline\delta^r (t)=\genfrac{}{}{1pt}{0}{r!}{\prod_{v\in\mathcal V(t)}r_v!}t_0,\]
where $t_0$ is the tree $t$ with all free edges removed.
\end{prop}
\begin{proof}
We have
\[\overline\delta^r  (t)=\sum_{(v_1,\ldots,v_r)\in\mathcal V(t)^r}\overline\delta_{v_1}\circ\cdots\circ\overline\delta_{v_r}(t).\]
A term in the right-hand side is equal to $t_0$ if and only if each vertex $v\in \mathcal V(t)$ appears exactly $r_v$ times in the tuple $(v_1,\ldots, v_r)$, otherwise the term is equal to zero. The number of such $r$-tuples is equal to the multinomial coefficient above.
\end{proof}
For any admissible cut $c$ of a tree $t\in\mop{PL}(A)$, let $\overline R_c(t)$ be the associated \textsl{full trunk}, which is given by the trunk $R^c(t)$ together with each cut edge replaced by a free edge. The weight of the full trunk is therefore equal to $r-1$, where $r$ is the number of edges belonging to the cut. In view of Proposition \ref{deltabar-max}, we have
\begin{equation}\label{trunk}
\overline\delta^r\big(\overline R^c(t)\big)=\|\overline t\|R^c(t),
\end{equation}
where $\overline t$ is a shorthand for $\overline R^c(t)$, and 
\begin{equation}\label{norm}
\|\overline t\|:=\genfrac{}{}{1pt}{0}{r!}{\prod_{v\in\mathcal V(\overline t)}r_v!}.
\end{equation}
\ignore{Let us define two coproducts $\overline\Delta_{\smop{BCK}}$ and $\overline\Delta_{\smop{LOT}}$ on $\overline {\mathcal H}_{\smop{BCK}}^A$ and $\overline {\mathcal H}_{\smop{LOT}}^A$  respectively, as follows. $\overline\Delta_{\smop{BCK}}$ is given by \eqref{coprod}, except that free edges incoming to each vertex remain, and that any edge relating $x\in V$ to $y\in W$ gives rise to a free edge attached to $y$ in the trunk $u\restr W$. The other coproduct is given by
\begin{equation}
\overline\Delta_{\smop{LOT}}(\mathbf x^{\mathbf k})=\sum_{r\ge 0}\,\sum_{\mathbf k=\mathbf k_1+\cdots +\mathbf k_r+\overline{\mathbf k}}\mathbf x_1^{\mathbf k_1}\odot\cdots\odot\mathbf x_r^{\mathbf k_r}\otimes \mathbf x^{\overline{\mathbf k}}
\end{equation}
where the inner sum runs over multisets $\{\mathbf k_1,\ldots,\mathbf k_r\}$ such that the remainder $\overline{\mathbf k}$ has only nonnegative components.\ignore{Let us extend derivation $\overline\partial$ to a derivation of $\overline {\mathcal H}_{\smop{LOT}}^A$ by making it a derivation with respect to the product $\odot$.} \begin{prop}
Both $\overline {\mathcal H}_{\smop{BCK}}^A$ and $\overline {\mathcal H}_{\smop{LOT}}^A$ are connected graded Hopf algebras, and both derivations $\overline\partial$ are also coderivations with respect to $\overline\Delta_{\smop{BCK}}$ and $\overline\Delta_{\smop{LOT}}$ respectively.
\end{prop}
\begin{proof}
Straightforward and left to the reader.
\end{proof}
}
\subsection{Graftings and cuts}
\begin{defn}\label{gft}
Let $r$ be a positive integer, let $F$ be an $A$-decorated rooted forest with $r$ connected components without free edges, and let $\overline t$ be an $A$-decorated rooted tree with $r_v$ free edges at vertex $v$, and a total number of $r$ free edges. A grafting of $F$ on $\overline t$ consists in grafting a choice of $r_v$ connected components of $F$ on the vertex $v$ for any $v\in\mathcal V(\overline t)$, and by removing the free edges.
\end{defn}
\noindent Denoting by $\mathcal G(F,\overline t)$ the set of graftings of $F$ on $\overline t$, we obviously have
\begin{equation}\label{cardinal-g}
\left\vert\mathcal G(F,\overline t)\right\vert = \|\overline t\|.
\end{equation}
For any vertex $x\in\mathcal V(\overline t)$ and for any grafting $b\in\mathcal G(F,\overline t)$, we'll denote by $F_b(x)$ the subforest of $F$ attached to $x$ via $b$.
\begin{prop}\label{sym-cut-graft}
Let $r$ be a positive integer, let $F$ and $\overline t$ as in Definition \ref{gft}, and let $t$ be an $A$-decorated rooted tree without free edges. Let $\mathcal G(t,F,\overline t)$ be the set of graftings of $F$ on $\overline t$ resulting in the tree $t$. Let $\mathcal C(t,F,\overline t)$ be the set of admissible cuts of $t$ such that $P^c(t)=F$ and $\overline R^c(u)=\overline t$. Then
\[\vert\mathcal C(t,F,\overline t)\vert=\frac{\sigma(t)}{\sigma(F)\sigma(\overline t)}\vert\mathcal G(t,F,\overline t)\vert.\]
\end{prop}
\begin{proof}
The group $\mop{Aut} t$ acts transitively on $\mathcal C(t,F,\overline t)$. The stabiliser of a cut $c$ will be denoted by $\mop{Aut}_c t$, and its cardinal by $\sigma_c(t)$. On the other hand, the group $\mop{Aut} \overline t\times\mop{Aut} F$ acts transitively on $\mathcal G(t, F,\overline t)$. To see this, consider two graftings $b,b'\in\mathcal G(t,F,\overline t)$: there exists a permutation $\alpha$ of $\mathcal V(\overline t)$ such that $F_b(x)$ and $F_{b'}\big(\alpha(x)\big)$ are isomorphic, and the permutation $\alpha$ necessarily comes from an automorphism of the tree $\overline t$. The stabiliser of $b$ is $\mop{Aut}_b\overline t\times\mop{Aut}_b F$, where $\mop{Aut}_b \overline t$ is the subgroup of those $\alpha\in\mop{Aut} \overline t$ such that $F_b\big(\alpha(x)\big)$ and $F_b(x)$ are isomorphic for any vertex $x$ of $\overline t$, and where $\mop{Aut}_b F=\prod_{x\in\mathcal V(\overline t)}\mop{Aut}F_b(x)$ is the subgroup of $\mop{Aut} F$ which respects the subforests $F_b(x)$. By the orbit-stabiliser theorem, we therefore have
\[\frac{\vert\mathcal C(t,F,\overline t)\vert}{\vert\mathcal G(t,F,\overline t)\vert}=\frac{\sigma(t)}{\vert\mop{Aut}_c(t)\vert}\frac{\vert\mop{Aut}_b(\overline t)\vert.\vert\mop{Aut}_b(F)\vert}{\sigma(F)\sigma(\overline t)}\]
We conclude by noticing the obvious isomorphism
\[\mop{Aut}_c(t)\sim \mop{Aut}_b(\overline t)\times\mop{Aut}_b(F).\]
\end{proof}

\section{Explicit description of the coproduct}\label{sect:main}
\subsection{A recursive formula for the coproduct}\label{par:rec}
Recall that the coproduct in $\mathcal H_{\smop{BCK}}^A$ admits a recursive definition, with respect to the degree, in terms of the cocycle operators $B_+^a$. This can be rewritten in terms of the operators $B^a_-$ as follows:
\begin{equation}\label{bckrec}
\Delta_{\smop{BCK}}(t)=\sum_{a\in A}(I\otimes B_+^a)\Delta_{\smop{BCK}}\big(B_-^a(t)\big)+t\otimes \un.
\end{equation}
We can now compute, using the Hopf algebra morphism property for $\jmath$:
\begin{eqnarray}\nonumber
(\jmath\otimes\jmath)\circ \Delta_{\smop{LOT}}(\mathbf x^{\mathbf k})&=&\Delta_{\smop{BCK}}\circ\jmath(\mathbf x^{\mathbf k})\\
&\mopl{=}_{\eqref{bckrec}}&\sum_{a\in A}(I\otimes B_+^a)\Delta_{\smop{BCK}}\big(B_-^a\circ\jmath(\mathbf x^{\mathbf k})\big)+\jmath(\mathbf x^{\mathbf k})\otimes \un\nonumber\\
&\mopl{=}_{\eqref{phibtr}}&\sum_{a\in A}(I\otimes B_+^a)\Delta_{\smop{BCK}}\circ\jmath\big(\overline L^a(\mathbf x^{\mathbf k})\big)+\jmath(\mathbf x^{\mathbf k})\otimes \un\nonumber\\
&=&\sum_{a\in A}(I\otimes B_+^a)\circ(\jmath\otimes\jmath)\circ \Delta_{\smop{LOT}}\big( \overline L^a(\mathbf x^{\mathbf k})\big)+\jmath(\mathbf x^{\mathbf k})\otimes \un\label{coprod-rec}.
\end{eqnarray}
\begin{remark}
The map $B_+^a\circ\jmath$ cannot be easily expressed like \eqref{phibtr} with $\jmath$ on the left, hence the above recursive expression cannot be further simplified.
\end{remark}
\subsection{An explicit formula for the coproduct}\label{par:main}
\begin{thm}\label{explicit}
The coproduct $\Delta_{\smop{LOT}}:\mathcal H_{\smop{LOT}}^A\to\mathcal H_{\smop{LOT}}^A\otimes\mathcal H_{\smop{LOT}}^A$ is the unique unital algebra morphism defined on the monomials by
\begin{equation}\label{coprod-N}
\Delta_{\smop{LOT}}(\mathbf x^{\mathbf k})=\sum_{r\ge 0}\,\sum_{\mathbf k=\mathbf k^1+\cdots +\mathbf k^r+\overline{\mathbf k},\, \smop{wt}(\mathbf k^j)=-1}\frac{\mathbf k!}{\sigma (\mathbf x^{\mathbf k^1}\odot\cdots\odot\mathbf x^{\mathbf k^r})\overline{\mathbf k}!}\,\mathbf x^{\mathbf k^1}\odot\cdots\odot\mathbf x^{\mathbf k^r}\otimes \overline{\partial}^r\mathbf x^{\overline{\mathbf k}},
\end{equation}
where the inner sum runs over multisets $\{\mathbf k^1,\ldots,\mathbf k^r\}$ of multi-indices such that $\mop{wt}(\mathbf k^j)=-1$ for any $j=1,\ldots,r$ and such that the remainder $\overline{\mathbf k}$, of weight $r-1$, has only nonnegative components $\overline k_j^a$.
\end{thm}
\begin{proof}
Let us first give an example: applying \eqref{coprod-N} to the monomial $x_{-1}^2x_0x_1$ yields
\begin{eqnarray*}
\Delta_{\smop{LOT}}(x_{-1}^2x_0x_1)&=&x_{-1}^2x_0x_1\otimes \un +\un\otimes x_{-1}^2x_0x_1+2x_{-1}\otimes \overline\partial(x_{-1}x_0x_1)
+2x_{-1}x_{0}\otimes \overline\partial(x_{-1}x_1)\\
&&+x_{-1}^2x_1\otimes\overline\partial x_0+x_{-1}\odot x_{-1}\otimes \overline\partial^2(x_0x_1)+2x_{-1}\odot x_{-1}x_0\otimes \overline\partial^2x_1\\
&=&x_{-1}^2x_0x_1\otimes \un +\un\otimes x_{-1}^2x_0x_1+2x_{-1}\otimes x_{-1}^2x_1+2x_{-1}\otimes x_{-1}x_0^2
+2x_{-1}x_{0}\otimes x_{-1}x_0\\
&&+x_{-1}^2x_1\otimes x_{-1}+3x_{-1}\odot x_{-1}\otimes x_{-1}x_0+2x_{-1}\odot x_{-1}x_0\otimes x_{-1}.
\end{eqnarray*}
We leave it to the reader to check that a repeated application of the recursive formula \eqref{coprod-rec}, starting from
\[\Delta_{\smop{LOT}}(x_{-1})=x_{-1}\otimes\un+\un\otimes x_{-1}\]
and computing successively the coproducts of $x_{-1}x_0$, $x_{-1}x_0^2$, $x_{-1}^2x_1$ and $x_{-1}^2x_0x_1$, gives the same result. The reader is invited to compare with Example 3.2 in \cite{BH24}. The differences between their coefficients and ours comes from the different convention for the symmetry factor.\\

\noindent Formula \eqref{coprod-N} for the coproduct admits the following alternative presentation:
\begin{equation}\label{coprod-N-bis}
\Delta_{\smop{LOT}}(\mathbf x^{\mathbf k})=\sum_{\mathbf c\in\smop{Adm}(\mathbf x^{\mathbf k})}P^{\mathbf c}(\mathbf x^{\mathbf k})\otimes R^{\mathbf c} (\mathbf x^{\mathbf k}),
\end{equation}
where $\mathbf c$ runs over all \textsl{admissible cuts} of the monomial $\mathbf x^{\mathbf k}$, i.e. all ways to cut $\mathbf k$ into $r+1$ multi-indices $\mathbf k^1,\ldots, \mathbf k^r,\overline{\mathbf k}$ for some integer $r\ge 0$, with $\mop{wt}(\mathbf k^j)=-1$ for any $j=1,\ldots,r$ and $\mop{wt}\overline{\mathbf k}=r-1$, so that
\[\mathbf x^{\mathbf k}= \mathbf x^{\mathbf k^1}\cdots \mathbf x^{\mathbf k^r}\mathbf x^{\overline{\mathbf k}}.\]
 In \eqref{coprod-N-bis}, we have set
\[P^{\mathbf c}(\mathbf x^{\mathbf k}):=\mathbf x^{\mathbf k^1}\odot\cdots\odot\mathbf x^{\mathbf k^r},\hskip 8mm R^{\mathbf c}(\mathbf x^{\mathbf k}):=\overline\partial^r\mathbf x^{\overline{\mathbf k}}.\]
In view of this, we shall denote by $\vert\mathbf c\vert$ the set of admissible  cuts $\mathbf c'$ such that $P^{\mathbf c'}(\mathbf x^{\mathbf k})=P^{\mathbf c}(\mathbf x^{\mathbf k})$ (and therefore $R^{\mathbf c'}(\mathbf x^{\mathbf k})=R^{\mathbf c}(\mathbf x^{\mathbf k})$), and by $\|\mathbf c\|$ the cardinal of this class. We clearly have
\[\|\mathbf c\|=\frac{\mathbf k!}{\overline{\mathbf k}!\sigma(\mathbf x^{\mathbf k^1}\odot\cdots\odot\mathbf x^{\mathbf k^r})}=\frac{\mathbf k!}{\mathbf k^1!\cdots\mathbf k^r!\overline{\mathbf k}!\,\sigma^{\smop{ext}}(\mathbf x^{\mathbf k^1}\odot\cdots\odot\mathbf x^{\mathbf k^r})},\]
whenever $P^{\mathbf c}(\mathbf x^{\mathbf k})=\mathbf x^{\mathbf k^1}\odot\cdots\odot\mathbf x^{\mathbf k^r}$, where $\sigma^{\smop{ext}}$ is the external symmetry factor (see Remark \ref{int-ext}). In analogy with trees, the \textsl{full trunk} will be defined by
\[\overline R^{\mathbf c}(\mathbf x^{\mathbf k}):=\mathbf x^{\overline{\mathbf k}}.\]
\begin{defn}\label{match}
Let $\mathbf c$ be an admissible cut of the monomial $\mathbf x^{\mathbf k}$, and let $c$ be an admissible cut of the decorated rooted tree $t$. We say that $c$ matches $\mathbf c$ and write $c\sim\mathbf c$ whenever
\begin{itemize}
\item $\mathbf x^{\mathbf k}=\Phi(t)$,
\item $P^{\mathbf c}(\mathbf x^{\mathbf k})=\Phi\big(P^c(t)\big)$.
\end{itemize}
\end{defn}
\noindent An admissible cut $c$ matches $\mathbf c$ if and only if it matches any element $\mathbf c'\in\vert\mathbf c\vert$. Note that the second condition implies $\overline R^{\mathbf c}(\mathbf x^{\mathbf k})=\overline\Phi\big(\overline R^c(t)\big)$. The theorem is an easy consequence of the following lemma:
\begin{lem}\label{lem:main}
For any monomial $\mathbf x^{\mathbf k}$ in $\mop N(A)$ and for any admissible cut $\mathbf c$ of  $\mathbf x^{\mathbf k}$, the following holds:
\begin{equation}\label{main-lemma}
\|\mathbf c\| (\jmath\otimes\jmath)\left(P^{\mathbf c}(\mathbf x^{\mathbf k})\otimes R^{\mathbf c}(\mathbf x^{\mathbf k})\right)=\sum_{t,\, \Phi(t)=\mathbf x^{\mathbf k}}\, \sum_{c\in\smop{Adm}(t),\, c\sim\mathbf c}
\frac{\sigma(\mathbf x^{\mathbf k})}{\sigma(t)}P^c(t)\otimes R^c(t).
\end{equation}
\end{lem}
\begin{proof}
Fix an admissible cut $\mathbf c$ of the monomial $\mathbf x^{\mathbf k}$. Denoting the left-hand side and the right-hand side of \eqref{main-lemma} by $\mathcal L$ and $\mathcal R$ respectively, we can compute:
\begin{eqnarray*}
\mathcal L &=&\frac{\mathbf k!}{\mathbf k^1!\cdots\mathbf k^r!\overline{\mathbf k}!\,\sigma^{\smop{ext}}(\mathbf x^{\mathbf k^1}\odot\cdots\odot\mathbf x^{\mathbf k^r})}\jmath(\mathbf x^{\mathbf k^1})\cdots \jmath(\mathbf x^{\mathbf k^r})\otimes \jmath\circ\overline\partial^r(\mathbf x^{\overline{\mathbf k}})\\
&=&(\mop{Id}\otimes\overline\delta^r)\left(\frac{\mathbf k!}{\mathbf k^1!\cdots\mathbf k^r!\overline{\mathbf k}!\,\sigma^{\smop{ext}}(\mathbf x^{\mathbf k^1}\odot\cdots\odot\mathbf x^{\mathbf k^r})}\jmath(\mathbf x^{\mathbf k^1})\cdots \jmath(\mathbf x^{\mathbf k^r})\otimes \jmath({\mathbf x}^{\overline{\mathbf k}})\right)\hbox{ (from Corollary \ref{deltabar-j})}\\
&=&\frac{\mathbf k!}{\mathbf k^1!\cdots\mathbf k^r!\overline{\mathbf k}!\,\sigma^{\smop{ext}}(\mathbf x^{\mathbf k^1}\odot\cdots\odot\mathbf x^{\mathbf k^r})}
(\mop{Id}\otimes\overline\delta^r)\left( \sum_{ \genfrac{}{}{0pt}{2}{(t^j)_{1,\ldots, r,}}  {\Phi(t^j)=\mathbf x^{\mathbf k^j} }}\ \sum_{\overline t,\, \Phi(\overline t)=\mathbf x^{\overline{\mathbf k}}}
\frac{\mathbf k^1!}{\sigma(t^1)}\cdots \frac{\mathbf k^r!}{\sigma(t^r)}\frac{\overline{\mathbf k}!}{\sigma(\overline t)}\, t^1\cdots t^r\otimes \overline t\right)\\
&=&\frac{\mathbf k!}{\sigma^{\smop{ext}}(\mathbf x^{\mathbf k^1}\odot\cdots\odot\mathbf x^{\mathbf k^r})}
(\mop{Id}\otimes\overline\delta^r)\left(\sum_{F,\, \Phi(F)=\mathbf x^{\mathbf k^1}\odot\cdots\odot\mathbf x^{\mathbf k^r}}\,\sum_{\overline t,\, \Phi(\overline t)=\mathbf x^{\overline{\mathbf k}}}\frac {\sigma^{\smop{ext}}(\mathbf x^{\mathbf k^1}\odot\cdots\odot\mathbf x^{\mathbf k^r})}{\sigma^{\smop{ext}}(F)}\frac{1}{\sigma^{\smop{int}}(F)\sigma(\overline t)}\,F\otimes \overline t\right)\\
&& \hbox{ (from Proposition \ref{prop:sym-forest})}\\
&=&\mathbf k!
(\mop{Id}\otimes\overline\delta^r)\left(\sum_{F,\, \Phi(F)=\mathbf x^{\mathbf k^1}\odot\cdots\odot\mathbf x^{\mathbf k^r}}\,\sum_{\overline t,\, \Phi(\overline t)=\mathbf x^{\overline{\mathbf k}}}\frac{1}{\sigma(F)\sigma(\overline t)}\,F\otimes \overline t\right).
\end{eqnarray*}
From \eqref{trunk}, \eqref{cardinal-g} and Proposition \ref{sym-cut-graft}, we therefore get
\begin{eqnarray*}
\mathcal L&=&\mathbf k!\sum_{F,\, \Phi(F)=\mathbf x^{\mathbf k^1}\odot\cdots\odot\mathbf x^{\mathbf k^r}}\,\sum_{\overline t,\, \Phi(\overline t)=\mathbf x^{\overline{\mathbf k}}}\|\overline t\| \frac{1}{\sigma(F)\sigma(\overline t)}\,F\otimes \overline t_0\\
&=& \mathbf k!\sum_{F,\, \Phi(F)=\mathbf x^{\mathbf k^1}\odot\cdots\odot\mathbf x^{\mathbf k^r}}\,\sum_{\overline t,\, \Phi(\overline t)=\mathbf x^{\overline{\mathbf k}}}\frac{\vert\mathcal G(F,\overline t)\vert}{\sigma(F)\sigma(\overline t)}\,F\otimes \overline t_0\\
&=& \mathbf k!\sum_{t,\,\Phi(t)=\mathbf x^{\mathbf k}}\,\sum_{F,\, \Phi(F)=\mathbf x^{\mathbf k^1}\odot\cdots\odot\mathbf x^{\mathbf k^r}}\,\sum_{\overline t,\, \Phi(\overline t)=\mathbf x^{\overline{\mathbf k}}}\frac{\vert\mathcal G(t,F,\overline t)\vert}{\sigma(F)\sigma(\overline t)}\,F\otimes \overline t_0\\
&=& \mathbf k!\sum_{t,\,\Phi(t)=\mathbf x^{\mathbf k}}\,\sum_{F,\, \Phi(F)=\mathbf x^{\mathbf k^1}\odot\cdots\odot\mathbf x^{\mathbf k^r}}\,\sum_{\overline t,\, \Phi(\overline t)=\mathbf x^{\overline{\mathbf k}}}\frac{\vert\mathcal C(t,F,\overline t)\vert}{\sigma(t)}\,F\otimes \overline t_0\\
&=& \sum_{t,\,\Phi(t)=\mathbf x^{\mathbf k}}\frac{\mathbf k!}{\sigma(t)}\sum_{c\in\smop{Adm}(t),\, c\sim\mathbf c}P^c(t)\otimes R^c(t)\\
&=&\mathcal R.
\end{eqnarray*}
\end{proof}
\noindent \textit{End of proof of Theorem \ref{explicit}.} By applying Lemma \ref{lem:main} and summing over all classes $\vert \mathbf c\vert$ of admissible cuts $\mathbf c$ of the monomial $\mathbf x^{\mathbf k}$,  the coproduct defined by \eqref{coprod-N}, or equivalently by \eqref{coprod-N-bis}, which we denote temporarily by $\Delta'_{\smop{LOT}}$, verifies
\[(\jmath\otimes\jmath)\circ \Delta'_{\smop{LOT}}=\Delta_{\smop{BCK}}\circ\jmath.\]
In view of the injectivity of $\jmath$, it therefore coincides with $\Delta_{\smop{LOT}}$, and the theorem is proven.
\end{proof}
\section{Extraction-contraction}
After recalling from \cite{CEM11} the extraction-contraction coproduct of (undecorated) rooted forests, we define in this section an extraction-contraction coproduct for multi-indices by Formula \eqref{ecbis} below, and we prove (Theorem \ref{ec-main}) that the injection $\jmath$ previously defined is a coalgebra morphism for both extraction-contraction coproducts. As a consequence $\jmath$ is a morphism of comodule-bialgebras, i.e. respects both extraction-contraction coproducts in addition to both Hopf algebra structures.
\subsection{Reminder on extraction-contraction of rooted forests}
Recall from \cite{CEM11} that a second coproduct $\Gamma$ on the algebra of (non-decorated) rooted forests makes it a commutative bialgebra, on which the Hopf algebra $\mathcal H_{\smop{BCK}}^{\{*\}}$ is a comodule-bialgebra. In particular, both coproducts are linked by a cointeraction diagram, dual to left distributivity. It is given by
\begin{equation}\label{ec}
\Gamma(t)=\sum_{s\subseteq t}s\otimes t/s,
\end{equation}
where the sum runs over the \textsl{covering subforests} of the forest $t$. A covering subforest is a partition of the set of vertices of $t$ into connected blocks, i.e. blocks in which any vertex can be reached from another by following edges of $t$. The notation $t/s$ stands for the corresponding contracted forest, obtained from $s$ by shrinking each block to a single vertex. A decorated version of this picture is available provided a commutative semigroup structure on the set of decorations is given, in order to decide how to decorate the vertices of the contracted forest.
\subsection{Free edges and extraction-contraction for multi-indices}
Formula \eqref{ec} can be precised as follows: the extraction of a covering subforest $s\subseteq t$ with $r(s)+1$ connected components gives rise to a forest $F(s)$ with $r(s)$ free edges, which are naturally in bijection with the edges of the contracted forest $t/s$. Note that two different covering subforests can give rise to the same forest after extraction. We denote by $F_0(s)$ the same forest without its free edges, given by
\[F_0(s)=\frac{1}{\|F(s)\|}\overline\delta^{r(s)}F(s),\]
where the coefficient $\|F(s)\|$ is given by the recipe \eqref{norm}. Formula \eqref{ec} therefore takes the following form:
\begin{equation}\label{ecter}
\Gamma(t)=\sum_{s\subseteq t}F_0(s)\otimes t/s.
\end{equation}

Such a coproduct $\Gamma$ does exist on the multi-index Hopf algebra \cite{L23}: we give in this section an explicit combinatorial formula for it, similar to \eqref{ec} the same way \eqref{coprod-N-bis} is similar to the Connes-Kreimer formula \eqref{bck-adm}. We stick to the nondecorated setting $A=\{*\}$ for simplicity, but a decorated version of this picture is also available when $A$ is endowed with a commutative semigroup structure. Details are left to the reader.\\

A \textsl{covering subforest} of a monomial $\mathbf x^{\mathbf k}=\prod_{j\ge -1}x_j^{k_j}$ (of weight $-1$) is a partition $\mathbf \pi$ of the corresponding multiset of variables into several multisets of total weight $\ge -1$, each of them giving rise to a monomial $\mathbf x^{\beta^j}$ for $j=1,\ldots, r(\pi)+1$. Here $r(\pi)+1$ is the number of multisets involved. We write $\pi\subseteq \mathbf x^{\mathbf k}$ for $\pi$ being a covering subforest of $\mathbf x^{\mathbf k}$. The corresponding monomial of monomials can be written as
\[\mathbb M(\mathbf \pi)=\mathbf x^{\mathbf \beta^1}\odot\cdots\odot \mathbf x^{\mathbf \beta^{r(\pi)+1}},\]
with $\mathbf \beta^1+\cdots +\mathbf \beta^{r(\pi)+1}=\mathbf k$. Its \textsl{reduced form} is defined by
\[\mathbb M_0(\mathbf \pi):=\frac{1}{\|\mathbf\pi\|}\overline\partial^{r(\pi)}\mathbb M(\mathbf \pi),\]
with
\begin{equation}\label{normofpi}
\|\mathbf\pi\|=\|\mathbb M(\mathbf\pi)\|:=\genfrac{}{}{0.5pt}{0}{r(\pi)!}{\prod_{j=1}^{r(\pi)+1} (\mop{wt}\mathbf \beta^j+1)!}.
\end{equation}
Our educated guess for the extraction-contraction coproduct on monomials is, on the model of \eqref{ecter}:
\begin{equation}\label{ecbis}
\Gamma(\mathbf x^{\mathbf k})=\sum_{\pi\subseteq \mathbf x^{\mathbf k}}\mathbb M_0(\pi)\otimes \mathbf x^{\mathbf k}/\pi,
\end{equation}
where the contracted monomial, of weight $-1$, is given by
\[\mathbf x^{\mathbf k}/\pi=\mathbf x^{\mathbf k}/\mathbb M(\pi):=\prod_{j=1}^{r(\pi)+1} x_{\smop{wt}\beta^j}.\]
\begin{thm}\label{ec-main}
The injection $\jmath$ is a coalgebra morphism with respect to both extraction-contraction coproducts \eqref{ec} and \eqref{ecbis}.
\end{thm}
\noindent The proof of Theorem \ref{ec-main} is postponed to Paragraph \ref{proof-ecmain} below.
\begin{cor}
The Hopf algebra $(\mathcal H_{\smop{LOT}}^{\{*\}},.,\Delta)$ is a comodule-bialgebra over the bialgebra $(\mathcal H_{\smop{LOT}}^{\{*\}},.,\Gamma)$.
\end{cor}
\subsection{Insertion and extraction of forests revisited}
\begin{defn}
Let $F$ be a rooted forest with $r(F)$ free edges and $r(F)+1$ connected components, and let $\overline t$ be a rooted tree with $r(F)+1$ vertices and $r(F)$ edges (hence without free edges). An insertion of $F$ inside $\overline t$ is a class of bijections
\[\tau:\{\hbox{edges of } \overline t\}\longrightarrow \{\hbox{free edges of } F\}\]
that induces a bijection $\widetilde\tau$ from $\mathcal V(\overline t)$ onto the set of connected components of $F$, in the sense that the restriction of $\tau$ to $E_v(\overline t)$ (the latter denoting the set of edges of $\overline t$ with bottom vertex $v$)  is the set of free edges of some connected component of $F$. Two such bijections $\tau$ and $\tau'$ are in the same class if and only if $\widetilde\tau=\widetilde{\tau'}$ and $\tau'=\varepsilon\circ\tau$, where $\varepsilon$ is a permutation of the free edges of $F$ that does not change their vertex.
\end{defn}
\begin{prop}\label{total-insertions}
Let $F$ be a rooted forest with $r(F)$ free edges and $r(F)+1$ connected components, let $\overline t$ be a rooted tree with $r(F)+1$ vertices and $r(F)$ edges, and let $\mathcal I(F,\overline t)$ be the set of insertions of $F$ inside $\overline t$. Let us denote by $e_v$ (resp. $r_w$) the number of edges of $\overline t$ with bottom vertex $v\in \mathcal V(\overline t)$ (resp. the number of free edges of $F$ attached to vertex $w\in\mathcal V(F)$). The following formula holds:
\[\vert \mathcal I(F,\overline t)\vert=\frac{\prod_{v\in\mathcal V(\overline t)}e_v!}{\prod_{w\in\mathcal V(F)}r_w!}\sigma\big(\Phi(\overline t)\big).\]
\end{prop}
\begin{proof}
The product in the numerator is the total number of bijections associated to a particular bijection from $\mathcal V(\overline t)$ onto the set of connected components of $F$. The denominator is the cardinal of each equivalence class, and the factor $\sigma\big(\Phi(\overline t)\big)$ is the number of permutations of $\mathcal V(\overline t)$ leaving invariant the fertility of each vertex. Proposition \ref{total-insertions} follows.
\end{proof}
\begin{prop}\label{evsi}
Let $F$ be a rooted forest with $r(F)$ free edges and $r(F)+1$ connected components, let $\overline t$ be a rooted tree with $r(F)+1$ vertices and $r(F)$ edges, and let $t$ be a rooted tree without free edges. Let $\mathcal I(t,F,\overline t)$ be the set of insertions of $F$ inside $\overline t$ such that the resulting tree is isomorphic to $t$, and let $\mathcal E(t,F,\overline t)$ be the set of extractions of $F$ from $t$ with contraction $\overline t$, i.e. covering subforests $s\subseteq t$ such that $F(s)\sim F$ and $t/s\sim \overline t$. Then we have
\[\frac{\mathcal E(t,F,\overline t)}{\mathcal I(t,F,\overline t)}=\frac{\sigma(t)}{\sigma(F)\sigma(\overline t)}.\]
\end{prop}
\begin{proof}
The group $\mop{Aut} t$ acts transitively on $\mathcal E(t,F,\overline t)$. On the other hand, the group $\mop{Aut} F\times\mop{Aut}\overline  t$ acts transitively on $\mathcal I(t,F,\overline t)$. Any insertion $\iota$ obviously gives rise to a covering subforest $s=s(\iota)$ of $t$. The stabilizer of $s\in\mathcal E(t,F,\overline t)$ is isomorphic to $\mop{Int}F(s)\times \mop{Aut}_s\overline t$, where $\mop{Int}F$ is the group of automorphisms of $F(s)$ leaving each connected component fixed, and where $\mop{Aut}_s\overline t$ is the group of automorphisms $\tau$ of $\overline t=t/s$ such that, for any $v\in\mathcal V(\overline t)$, both connected components of $s$ corresponding to $v$ and $\tau(v)$ are isomorphic. This is precisely the stabilizer of $\iota$. By the orbit-stabilizer theorem we therefore get
\begin{eqnarray*}
\frac{\mathcal E(t,F,\overline t)}{\mathcal I(t,F,\overline t)}&=&\frac{\sigma(t)/\vert\mop{Stab} s\vert}{\sigma(F)\sigma(\overline t)/\vert\mop{Stab} \iota\vert}\\
&=&\frac{\sigma(t)}{\sigma(F)\sigma(\overline t)}.
\end{eqnarray*}
\end{proof}
\subsection{Proof of Theorem \ref{ec-main}}\label{proof-ecmain}
We have
\begin{eqnarray*}
\mathcal R:=\Gamma\big( \jmath(\mathbf x^{\mathbf k})\big)&=& \sum_{t,\,\Phi(t)=\mathbf x^{\mathbf k}}\frac {k!}{\sigma(t)}\Gamma(t)\\
&=&\sum_{t,\,\Phi(t)=\mathbf x^{\mathbf k}}\frac {k!}{\sigma(t)}\sum_{s\subseteq t}\frac{1}{\|F(s)\|}\overline\delta^{r(s)}F(s)\otimes t/F(s).
\end{eqnarray*}
On the other hand,
\begin{eqnarray*}
\mathcal L:=(\jmath\otimes\jmath)\big(\Gamma(\mathbf x^{\mathbf k})\big)
&=&\sum_{\pi\subseteq\mathbf x^{\mathbf k}}\frac {1}{\|\pi\|}\jmath\big(\overline\partial^{r(\pi)}\mathbb M(\pi)\big)\otimes \jmath(\mathbf x^{\mathbf k}/\pi)\\
&=&\sum_{\mathbb M}\,\sum_{\pi\subseteq\mathbf x^{\mathbf k},\, \mathbb M(\pi)=\mathbb M}\frac {1}{\|\pi\|}\jmath\big(\overline\partial^{r(\pi)}\mathbb M\big)\otimes \jmath(\mathbf x^{\mathbf k}/\mathbb M).
\end{eqnarray*}
Here, the external sum runs over monomials of monomials $\mathbb M=\mathbf x^{\beta^1}\odot\cdots\odot\mathbf x^{\beta^{r(\mathbb M(\pi))+1}}$, with $\mop{wt}\beta^j\ge -1$ and $\beta^1+\cdots +\beta^{r(\mathbb M)+1}=\mathbf k$. Using $r(\mathbb M)=r(\pi)$ whenever $\mathbb M=\mathbb M(\pi)$ and following the same lines than in the proof of Lemma \ref{lem:main}, we compute:
\begin{eqnarray*}
\mathcal L&=&\sum_{\mathbb M} \frac{\mathbf k!}{\sigma(\mathbb M)}\frac {1}{\|\mathbb M\|}\jmath\big(\overline\partial^{r(\mathbb M)}\mathbb M\big)\otimes \jmath(\mathbf x^{\mathbf k}/\mathbb M)\\
&=&\sum_{\mathbb M} \frac{\mathbf k!}{\sigma(\mathbb M)}\frac {1}{\|\mathbb M\|}\overline\delta^{r(\mathbb M)}\jmath(\mathbb M)\otimes \jmath(\mathbf x^{\mathbf k}/\mathbb M)\\
&=&\sum_{\mathbb M} \frac{\mathbf k!}{\sigma(\mathbb M)}\frac {1}{\|\mathbb M\|}\overline\delta^{r(\mathbb M)}\left(\jmath(\mathbf x^{\beta^1})\cdots\jmath(\mathbf x^{\beta^{r(\mathbb M)+1)}})\right)\otimes \jmath(\mathbf x^{\mathbf k}/\mathbb M)\\
&=&\sum_{\mathbb M} \frac{\mathbf k!}{\sigma(\mathbb M)\|\mathbb M\|}\overline\delta^{r(\mathbb M)}
\sum_{t^1,\ldots, t^{r(\mathbb M)+1},\overline t,\, \Phi(t^j)=\mathbf x^{\beta^j}, \, \Phi(\overline t)=\mathbf x^{\mathbf k}/\mathbb M}\frac{\beta^1!\cdots \beta^{r(\mathbb M)+1}!\sigma(\mathbf x^{\mathbf k}/\mathbb M)}{\sigma(t^1)\cdots\sigma(t^{r(\mathbb M)+1})\sigma(\overline t)}\,t^1\cdots t^{r(\mathbb M)+1}\otimes \overline t\\
&=&\sum_{\mathbb M} \frac{\mathbf k!}{\sigma(\mathbb M)\|\mathbb M\|}\overline\delta^{r(\mathbb M)}
\sum_{F,\overline t,\, \Phi(F)=\mathbb M, \, \Phi(\overline t)=\mathbf x^{\mathbf k}/\mathbb M}\frac{\beta^1!\cdots \beta^{r(\mathbb M)+1}!\sigma(\mathbf x^{\mathbf k}/\mathbb M)}{\sigma(t^1)\cdots\sigma(t^{r(\mathbb M)+1})\sigma(\overline t)}\frac{\sigma^{\smop{ext}(\mathbb M)}}{\sigma^{\smop{ext}(F)}}\, F\otimes \overline t\\
&=&\sum_{\mathbb M} \frac{\mathbf k!}{\|\mathbb M\|}\overline\delta^{r(\mathbb M)}
\sum_{F,\overline t,\, \Phi(F)=\mathbb M, \, \Phi(\overline t)=\mathbf x^{\mathbf k}/\mathbb M}\frac{\sigma(\mathbf x^{\mathbf k}/\mathbb M)}{\sigma(F)\sigma(\overline t)}\, F\otimes \overline t.
\end{eqnarray*}
From Proposition \ref{total-insertions} and Proposition \ref{evsi} we therefore get
\begin{eqnarray*}
\mathcal L&=&\sum_{\mathbb M} \frac{\mathbf k!}{\|\mathbb M\|}\overline\delta^{r(\mathbb M)}
\sum_{F,\overline t,\, \Phi(F)=\mathbb M, \, \Phi(\overline t)=\mathbf x^{\mathbf k}/\mathbb M}\mathcal I(F,\overline t)\frac{\prod_{w\in\mathcal V(F)}r_w!}{\left(\prod_{v\in\mathcal V(\overline t)}e_v!\right)\sigma\big(\Phi(\overline t)\big)}\frac{\sigma(\mathbf x^{\mathbf k}/\mathbb M)}{\sigma(F)\sigma(\overline t)}\, F\otimes \overline t\\
&=&\sum_{t,\, \Phi(t)=\mathbf x^{\mathbf k}}\,\sum_{\mathbb M} \frac{\mathbf k!}{\|\mathbb M\|}\overline\delta^{r(\mathbb M)}
\sum_{F,\overline t,\, \Phi(F)=\mathbb M, \, \Phi(\overline t)=\mathbf x^{\mathbf k}/\mathbb M}\mathcal I(t,F,\overline t)\frac{\prod_{w\in\mathcal V(F)}r_w!}{\left(\prod_{v\in\mathcal V(\overline t)}e_v!\right)\sigma\big(\Phi(\overline t)\big)}\frac{\sigma(\mathbf x^{\mathbf k}/\mathbb M)}{\sigma(F)\sigma(\overline t)}\, F\otimes \overline t\\
&=&\sum_{t,\, \Phi(t)=\mathbf x^{\mathbf k}}\,\sum_{\mathbb M} \frac{\mathbf k!}{\|\mathbb M\|}\overline\delta^{r(\mathbb M)}
\sum_{F,\overline t,\, \Phi(F)=\mathbb M, \, \Phi(\overline t)=\mathbf x^{\mathbf k}/\mathbb M}\mathcal E(t,F,\overline t)\frac{\prod_{w\in\mathcal V(F)}r_w!}{\prod_{v\in\mathcal V(\overline t)}e_v!}\frac{1}{\sigma(F)\sigma(\overline t)}
\frac{\sigma(F)\sigma(\overline t)}{\sigma(t)}\, F\otimes \overline t\\
&=&\sum_{t,\, \Phi(t)=\mathbf x^{\mathbf k}} \frac{\mathbf k!}{\sigma(t)}\,\sum_{\mathbb M} \frac{1}{\|\mathbb M\|}\overline\delta^{r(\mathbb M)}
\sum_{F,\overline t,\, \Phi(F)=\mathbb M, \, \Phi(\overline t)=\mathbf x^{\mathbf k}/\mathbb M}\mathcal E(t,F,\overline t)\frac{\prod_{w\in\mathcal V(F)}r_w!}{\prod_{v\in\mathcal V(\overline t)}e_v!}\, F\otimes \overline t.
\end{eqnarray*}
From \eqref{norm} and \eqref{normofpi}, we finally get
\begin{eqnarray*}
\mathcal L&=&\hskip -5mm \sum_{t,\, \Phi(t)=\mathbf x^{\mathbf k}} \frac{\mathbf k!}{\sigma(t)} \sum_{\mathbb M} \frac{\prod_{j=1}^{r(\mathbb M)+1} (\mop{wt}\beta^j+1)!}{r(\mathbb M)!}\,\overline\delta^{r(\mathbb M)}
\hskip -10mm\sum_{F,\overline t,\, \Phi(F)=\mathbb M, \, \Phi(\overline t)=\mathbf x^{\mathbf k}/\mathbb M}\frac{r(F)!}{\prod_{w\in\mathcal V(F)}r_w!}\frac{1}{\|F\|}\mathcal E(t,F,\overline t)\frac{\prod_{w\in\mathcal V(F)}r_w!}{\prod_{v\in\mathcal V(\overline t)}e_v!}\, F\otimes \overline t\\
&=& \sum_{t,\, \Phi(t)=\mathbf x^{\mathbf k}} \frac{\mathbf k!}{\sigma(t)}\,\sum_{\mathbb M}\, \overline\delta^{r(\mathbb M)}
\sum_{F,\overline t,\, \Phi(F)=\mathbb M, \, \Phi(\overline t)=\mathbf x^{\mathbf k}/\mathbb M}\frac{1}{\|F\|}\mathcal E(t,F,\overline t)\, F\otimes \overline t\\
&=&\mathcal R,
\end{eqnarray*}
which proves Theorem \ref{ec-main}.

\end{document}